\documentclass{lmcs}
\pdfoutput=1

\usepackage{lastpage}
\renewcommand{\dOi}{14(3:3)2018}
\lmcsheading{}{1--\pageref{LastPage}}{}{}%
{Apr.~24,~2017}{Jul.~20,~2018}{}

\usepackage{hyperref}
\hypersetup{hidelinks}

\theoremstyle{plain}

\usepackage{yfonts,array}
\usepackage{amsmath,amssymb,amsthm,units,stmaryrd}
\usepackage{upgreek,xcolor,bm,tikz,stackrel}
\usepackage[inline]{enumitem}
\usepackage{graphicx}
\usepackage[all]{xy}
\usepackage{tikz-cd}
\usetikzlibrary{arrows}

\newcommand{\cl }{\mathcal }
\newcommand{\fr }{\mathfrak }

\def\nb{\blacksquare}
\def\pb{\blacklozenge}

\newcommand{\lm}{\vec{\mathfrak Q}}
\newcommand{\domlm}{|\vec {\mathfrak Q}|}
\newcommand{\Tlm}{\vec {\mathcal T}_\peq}
\newcommand{\Slm}{\vec {S}}

\newcommand{\iiff}{\mathop\Leftrightarrow}
\newcommand{\ineg}{{\sim}}
\newcommand{\imp}{\mathop \Rightarrow}
\newcommand{\minsym}{{\sf m}}

\newcommand{\DTLm}{{\sf DTL}^\minsym }
\newcommand{\ITLm}{{\sf ITL}^\minsym}

\newcommand{\itlaipl}{{\sf ITL}^{\sf p} }
\newcommand{\itlepl}{{\sf ITL}^{\sf e} }

\newcommand{\itle}{{\sf ITL}^{\sf e}_\ps }
\newcommand{\itlk}{{\sf ITL}^{\sf c}_\nec}
\newcommand{\icltl}{{\sf ITL}^{\sf c}_\ps }
\newcommand{\itlcplus}{{\sf ITL}^{\sf c}}

\newcommand{\lanclass}{{\mathcal L}^C}
\newcommand{\lanfull}{{\mathcal L}^{I}}
\newcommand{\lank}{{\mathcal L}^I_\nec}
\newcommand{\landi}{{\mathcal L}^I_\ps }

\newcommand{\lgt}[1]{\|#1\|}
\newcommand{\nx}{{\circ}}
\newcommand{\nec}{\Box}
\newcommand{\ps}{\Diamond}
\newcommand{\val}[1]{\lb #1 \rb}

\newcommand{\cval}[1]{\lb #1 \rb^C}
\newcommand{\moments}[1]{\mathbb M_{ #1 }}
\newcommand{\irr}[1]{\mathbb I_{ #1 }}
\newcommand{\type}[1]{\mathbb T_{ #1 }}
\newcommand{\peq}{\preccurlyeq}
\def \bangle{ \atopwithdelims [ ]}
\newcommand{\add}[2]{{{#1}\bangle {#2}}}
\newcommand{\circrel}{\mathrel S_\Sigma}
\newcommand{\basic}[2]{\mathord\downarrow_{#1}#2}
\newcommand{\acc}{\peq}
\newcommand{\cca}{\succcurlyeq}
\newcommand{\defect}[1]{\partial {#1}}
\newcommand{\Root}[1]{r_{#1}}
\newcommand{\dom}[1]{{\rm dom}({#1})}
\newcommand{\fw}{{\mathfrak w}}
\newcommand{\fv}{{\mathfrak v}}
\newcommand{\fu}{{\mathfrak u}}
\newcommand{\redu}{\mathrel\unlhd}
\def\lb{\left\llbracket}
\def\rb{\right\rrbracket}
\def\<{\left (}

\def\>{\right )}
\def\({\left (}
\def\){\right )}

\def\cbra{\left \{}
\def\cket{\right \}}

\DeclareSymbolFont{AMSb}{U}{msb}{m}{n}
\DeclareMathSymbol{\N}{\mathbin}{AMSb}{"4E}
\DeclareMathSymbol{\Z}{\mathbin}{AMSb}{"5A}
\DeclareMathSymbol{\R}{\mathbin}{AMSb}{"52}
\DeclareMathSymbol{\Q}{\mathbin}{AMSb}{"51}
\DeclareMathSymbol{\I}{\mathbin}{AMSb}{"49}
\DeclareMathSymbol{\C}{\mathbin}{AMSb}{"43}

\newtheorem{defn}[thm]{Definition}
\newtheorem{corollary}[thm]{Corollary}
\newtheorem{exam}[thm]{Example}
\newtheorem{remark}[thm]{Remark}
\newtheorem{question}[thm]{Question}

\begin{document}

\title{The intuitionistic temporal logic of dynamical systems}

\author{David Fern\'andez-Duque}	
\address{Department of Mathematics, Ghent University, Belgium}	
\email{David.FernandezDuque@UGent.be}  





\begin{abstract}
A {\em dynamical system} is a pair $(X,f)$, where $X$ is a topological space and $f\colon X\to X$ is continuous. Kremer observed that the language of propositional linear temporal logic can be interpreted over the class of dynamical systems, giving rise to a natural intuitionistic temporal logic. We introduce a variant of Kremer's logic, which we denote $\icltl$, and show that it is decidable. We also show that minimality and Poincar\'e recurrence are both expressible in the language of $\icltl$, thus providing a decidable logic capable of reasoning about non-trivial asymptotic behavior in dynamical systems.
\end{abstract}

\maketitle

\section{Introduction}

{\em Dynamical (topological) systems} are mathematical models of change or movement over time. Formally, they are structures of the form $\mathfrak X=(X,\mathcal T,f)$, where $(X,\mathcal T)$ is a topological space and $f\colon X\to X$ is a continuous function \cite{akin}. This rather broad definition allows them to be found in multiple disciplines, including mathematics, physics and biology. 
The theory of such systems can be readily formalized in standard mathematical foundations such as set theory or even second-order arithmetic, but these are undecidable and subject to G\"odel's incompleteness theorems \cite{Godel1931}. On the other hand, G\"odel's proof relies on number-theoretic considerations which are not necessarily relevant to the study of dynamical systems.
In view of this, the {\em dynamic topological logic} project aimed to develop a logical framework for formal reasoning about topological dynamics which was decidable and complete yet powerful enough to prove non-trivial theorems. Unfortunately, the project was set back by a series of negative results, including the undecidability of full dynamic topological logic ($\sf DTL$), leading to a search for decidable variants of $\sf DTL$ which remained expressive enough for its intended applications.
Our goal in this paper is to present such a variant.  Following a suggestion of Kremer, we recast dynamic topological logic as an intuitionistic temporal logic, and show it to be a decidable logic powerful enough to reason about non-trivial recurrence phenomena in arbitrary dynamical systems.

\subsection{Dynamic topological logic}

As observed by Artemov et al.~\cite{arte}, it is possible to reason about dynamical systems within modal logic; the topological structure can be represented via a modality, which we will denote $\blacksquare$, interpreted as an interior operator in the sense of Tarski \cite{tarski}, while the action of $f$ can be captured by a temporal-like modality, which we will denote $\nx$. They introduced the logic $\sf S4C$, and proved that it is decidable, as well as sound and complete for the class of all dynamical systems. Kremer and Mints \cite{kmints} considered a similar logic, called $\sf S4H$, and also showed it to be sound and complete for the class of dynamical systems where $f$ is a homeomorphism.

The latter also suggested adding the `henceforh' operator, $\nec$, from linear temporal logic ($\sf LTL$) \cite{temporal}. This allowed them to represent asymptotic phenomena, including the non-trivial Poincar\'e recurrence theorem, which we discuss in \S \ref{SecProb}. The resulting tri-modal system was called {\em dynamic topological logic} ($\sf DTL$), and the positive results concerning $\sf S4C$ and $\sf S4H$ led to the expectation that $\sf DTL$ may also be well-behaved, possibly leading to applications in e.g.~specialized computer-assisted proof. Unfortunately, it was soon shown by Konev et al.~that dynamic topological logic is undecidable over the class of all dynamical systems \cite{konev}. Although the axiomatization of $\sf DTL$ proposed in \cite{kmints} is incomplete \cite{DFDTOCL}, we later provided a sound and complete axiomatization \cite{dtlaxiom}. However, Konev et al.~showed that $\sf DTL$ over the class of dynamical systems with a homeomorphism is not even computably enumerable \cite{wolter}.

This led to a search for decidable variants of $\sf DTL$ which retained the capacity for reasoning about the asymptotic behavior of dynamical systems. Gabelaia et al.~\cite{pml} proposed dynamic topological logics with finite, but unbounded, time and showed them to be decidable, although not in primitive recursive time. Kremer instead proposed a restriction to dynamical systems where the topology is a partition \cite{s5}, while we considered interpretations over minimal systems, which we will discuss in \S \ref{SecMinimal}; the $\sf DTL$'s obtained in the latter two cases are decidable. 

However, as a general rule, all decidable variants of $\sf DTL$ with continuous functions that are currently known are obtained by either restricting the class of dynamical systems over which they are interpreted, or restricting the logics to reason about finitely many iterations of $f$.\footnote{An exception to this is the `temporal over topological' fragment discussed in \cite{kmints}, but this fragment is too weak to express continuity.} This makes them unsuitable to capture asymptotic behavior of arbitrary dynamical systems, which motivated interest in $\sf DTL$. 

\subsection{Intuitionistic dynamic topological logic}

There is another variant of $\sf DTL$ which does not have either of these restrictions, yet whose decidability was never settled: Kremer's intuitionistic version of $\sf DTL$, proposed in unpublished work \cite{KremerIntuitionistic}, and here denoted $\itlk$. It is well-known that propositional intuitionistic logic can be seen as a fragment of $\sf S4$ via the G\"odel-Tarski translation \cite{tarski} (see \S \ref{SecComp} for details), and indeed they share very similar semantics. In particular, intuitionistic logic can be interpreted topologically. One can use this idea to present a version of dynamic topological logic which removes the modality $\blacksquare$, and instead interprets implication intuitionistically. A feature of such semantics is that all formulas are interpreted by {\em open} sets, and as Kremer observed, the truth condition for $\nx\varphi$ will preserve the openness of the truth valuation when $f$ is a continuous function. On the other hand, Kremer also observed that the classical truth condition for $\nec\varphi$ did not always produce open sets, and thus he proposed to take the interior of this classical interpretation (see \S \ref{SecBasic} for details). 

We will follow Kremer in using intuitionistic temporal logic to reason about dynamical systems, but we will replace $\nec$ by `eventually', denoted $\ps$; note that the two are not inter-definable intuitionistically \cite{Balbiani2017}. Working with $\ps$ is convenient because the classical interpretation of $\ps\varphi$ will indeed yield open sets. We will also add a universal modality, and call the resulting logic $\icltl$ {\em (Intuitionistic Temporal Logic of Continuous functions)}.

\subsection{Intuitionistic temporal logics}

The logic $\itlk$ is not the first intuitionistic variant of temporal logic. Ewald considered intuitionistic logic with `past' and `future' tenses \cite{Ewald} and Davies suggested an intuitionistic temporal logic with $\nx$ \cite{Davies96}, which was endowed with Kripke semantics and a complete deductive system by Kojima and Igarashi \cite{KojimaNext}. Logics with $\nx,\nec$ over bounded time were later studied by Kamide and Wansing \cite{KamideBounded}, and over infinite time by Nishimura \cite{NishimuraConstructivePDL}, who provided a sound and complete axiomatization for an intuitionistic variant of propositional dynamic logic $\sf PDL$. For logics with $\nx,\ps,\nec$, Balbiani and Di\'eguez~\cite{BalbianiDieguezJelia} axiomatized the `here-and-there' variant of $\sf LTL$, and with Boudou and Di\'eguez we showed that intuitionistic $\sf LTL$ over expanding frames, denoted $\itlepl$, is decidable \cite{Boudou2017}. 
Each of these logics use semantics based on bi-relational models for intuitionistic modal logic, studied systematically by Simpson \cite{Simpson94}. Topological semantics for modal logic in general, and for temporal logic with `past' and `future' in particular, have also been studied by Davoren et al.~\cite{Davoren09,DavorenIntuitionistic}.

\subsection{Our results} The main contribution in this article is that $\icltl$ is decidable, which we prove using techniques first introduced in \cite{me2}. There are a few simplifications due to the intuitionistic semantics, which allows us to work only with open sets in topological models, as well as partial orders rather than preorders in relational models. 
The semantics of $\sf DTL$ is reduced to countable (non-deterministic) quasimodels in~\cite{me2}, but those of $\icltl$ can be reduced to {\em finite} quasimodels; this is, of course, good news, but requires new combinatorial work with what we call irreducible moments. We will also show that minimality, as well as the Poincar\'e recurrence theorem, can be represented in $\icltl$, thus providing a promising new direction for the dynamic topological logic project. Nevertheless, this is an exploratory work, and many crucial questions remain open: our techniques do not provide us with an axiomatization, and the current decision procedure is super-exponential, meaning that further advances are needed to obtain a tractable logic.

\subsection*{Layout} The layout of the article is as follows.
\begin{itemize}[align=left]

\item[\S \ref{SecTopre}.] A brief review of some basic notions, including topology and its relation to partial orders.%

\item[\S \ref{SecBasic}.] Introduction to the logic $\icltl$ and its intended semantics using (dynamical topological) models.%

\item[\S \ref{SecComp}.] A comparison of our logic $\icltl$ to some related logics found in the literature.%

\item[\S \ref{SecNDQ}.] Introduction to labeled spaces and systems. These generalize both models and quasimodels, and the latter are defined as a special case.%

\item[\S \ref{SecGen}.] Construction of the limit model of a quasimodel, used to prove that $\icltl$ is sound for the class of quasimodels.%

\item[\S \ref{SecFour}.] Introduction to moments, which are essentially the points in the `weak canonical quasimodel' $\moments\Sigma$.%

\item[\S \ref{SecSim}.] Discussion of simulations, the essential tool for extracting quasimodels from models.%

\item[\S \ref{SecIrred}.] Introduction to the structure $\irr\Sigma$ of irreducible moments as a finite substructure of $\moments\Sigma$.

\item[\S \ref{SecSDS}.] Proof that a formula is falsifiable if and only if it is falsifiable over an effectively bounded quasimodel, obtaining the decidability of $\icltl$ as a corollary.

\item[\S \ref{SecSpec}.] Discussion on the capacity of $\icltl$ to reason about minimality and Poincar\'e recurrence, including a proof that the full $\itlcplus$ is decidable on the class of minimal systems.

\item[\S \ref{SecConc}.] Concluding remarks, including several open questions.

\end{itemize}

\section{Binary Relations and Topology}\label{SecTopre}

In this section we establish some of the notation we will use and recall some basic notions from topology. In particular, we discuss Alexandroff spaces, which link partial orders and topological spaces.

\subsection{Sets and relations}

We follow fairly standard conventions for sets and relations, which we briefly outline for reference.

The first infinite ordinal will be denoted $\omega$, the cardinality of a set $A$ will be denoted $\# A$, and its powerset will be denoted $\mathcal P(A)$. For a binary relation $R\subseteq A\times B$ and $X\subseteq A$, we write $R(X)$, or simply $R X$, for $\{y\in B: \exists x\in X \ x\mathrel R y\}$, and we write $R^{-1}$ for the inverse of $R$, that is, $\{(y,x)\in B\times A : (x,y)\in R\}$. We write $\dom R$ for $R^{-1}(B)$ and ${\rm rng}(R)$ for $R(A)$. $R$ is {\em total} if $\dom R=A$, and {\em surjective} if ${\rm rng}(R)=B$.

For any $X\subseteq A$, we define the restriction of $R$ to $X$ by $R\upharpoonright X =R\cap (X\times B)$, except when $A=B$, in which case $R\upharpoonright X =R\cap ( X\times X)$. If, moreover, $S\subseteq B\times C$, the composition of $R$ and $S$ is $S\circ R\subseteq A\times C$, which we may sometimes denote $SR$. If $B=A$ (so that $R\subseteq A\times A$), recall that $R$ is
\begin{enumerate}[label=(\alph*)]

\item {\em reflexive,} if for every $x\in A$, $x\mathrel R x$;

\item {\em antisymmetric,} if for every $x,y\in A$, if $x\mathrel R y$ and $y\mathrel R x$, then $x=y$;

\item {\em transitive,} if for every $x,y\in A$, if $x\mathrel R y$ and $y\mathrel R z$, then $x\mathrel R z$;

\item {\em linear,} if, for all $x,y\in A$, either $x \mathrel R y$ or $y \mathrel R x$, and

\item {\em serial,} if for every $x\in A$ there is $y\in A$ such that $x\mathrel R y$.

\end{enumerate}
A transitive, reflexive relation is a {\em preorder} and a transitive, reflexive, antisymmetric relation is a {\em partial order.}

A {\em poset} is a pair $\mathfrak A=(|\mathfrak A|,\peq_\mathfrak A)$, where $|\mathfrak A|$ is any set and ${\peq}_\mathfrak A \subseteq |\mathfrak A|\times |\mathfrak A|$ is a partial order; we will generally adopt the convention of denoting the domain of a structure $\mathfrak A$ by $|\mathfrak A|$. If $X\subseteq |\fr A|$, $\fr A\upharpoonright X$ is the structure obtained by restricting each component of $\fr A$ to $X$. Define $\mathop\downarrow_\mathfrak A a=\cbra b:b\peq_\mathfrak A a\cket$ and $\mathop\uparrow_\mathfrak A a=\cbra b:b\cca_\mathfrak A a\cket$. We will generally write $\peq,\downarrow,\uparrow$ instead of $\peq_\mathfrak A,\downarrow_\mathfrak A,\uparrow_\mathfrak A$ when this does not lead to confusion. We will also use the notation $a\prec b$ for $a\peq b$ but $b\not\peq a$, while $a\prec^1 b$ means that $a\prec b$ and there is no $c\in |\fr A|$ such that $a\prec c \prec b$. An element $a\in |\fr A|$ is {\em greatest} if $b\peq a$ for all $b\in |\fr A|$, and {\em maximal} if $b\not\succ a$ for any $b\in |\fr A|$. {\em Least} and {\em minimal} are defined analogously.

A partial order $\mathfrak A$ is a {\em tree} if it has a greatest element $r$, and for every $x\in |\mathfrak A|$, $\mathop\uparrow x$ is finite and linearly ordered. Note that in our presentation, the leaves of a tree are its minimal elements.

\subsection{Notions from topology}

Here we recall the basic definitions from topology we will use, including Alexandroff spaces, which link topological spaces and partial orders. A more in-depth introduction can be found in a standard text such as \cite{Munkres}.

\begin{defn}
A {\em topological space} is a pair $\mathfrak X=\<|\mathfrak X|,\mathcal{T}_\mathfrak X\>,$ where $|\mathfrak X|$ is a set and $\mathcal T_\mathfrak X$ a family of subsets of $|\mathfrak X|$ satisfying
\begin{enumerate}
\item $\varnothing,|\mathfrak X|\in \mathcal T_\mathfrak X$;
\item if $U,V\in \mathcal T_\mathfrak X$ then $U\cap V\in \mathcal T_\mathfrak X$ and
\item if $\mathcal O\subseteq\mathcal T_\mathfrak X$ then $\bigcup\mathcal O\in\mathcal T_\mathfrak X$.
\end{enumerate}
The elements of $\mathcal T_\mathfrak X$ are called {\em open sets}. Complements of open sets are {\em closed sets}. If $x\in U\in \cl T_\fr X$, we say that $U$ is a {\em neighborhood} of $x$.

If $\fr X,\fr Y$ are topological spaces, a function $f\colon |\mathfrak X| \to |\mathfrak Y|$ is {\em continuous} if $f^{-1}(V)$ is open whenever $V\subseteq |\mathfrak Y|$ is open, and {\em open} if $f(U)$ is open whenever $U\subseteq |\fr X|$ is open. A {\em homeomorphism} is an open, continuous bijection.
\end{defn}

It is readily verified that the composition of continuous functions is continuous, and similarly the composition of open functions is open; this will be useful throughout the text.
If $\mathfrak X$ is a topological space and $A\subseteq |\mathfrak X|$ is not open, $A$ can still be approximated by an open set from below.
To be precise, given a set $A\subseteq |\mathfrak X|$, its {\em interior}, denoted $A^\circ$, is defined by
\[A^\circ=\bigcup\cbra U\in\mathcal T_\mathfrak X:U\subseteq A\cket.\]
Similarly, we define the closure $\overline A$ as $|\mathfrak X|\setminus(|\mathfrak X|\setminus A)^\circ$; this is the smallest closed set containing $A$.

Perhaps the most standard example of a topological space is given by the real line $\mathbb R$, where $U\subseteq \mathbb R$ is open if, whenever $x\in U$, there is $\varepsilon > 0$ such that $(x-\varepsilon,x+\varepsilon ) \subseteq U$. More generally, let $(X,d)$ be a {\em metric space;} that is, $d\colon X\times X\to [0,\infty)$ satisfies 
\begin{enumerate}[label=(\alph*)]

\item $d(x,y)=0$ if and only if $x=y$,

\item  $d(x,y)=d(y,x)$, and

\item  the triangle inequality, $d(x,z)\leq d(x,y)+d(y,z)$.

\end{enumerate}
Then, for $\varepsilon>0$ and $x\in X$, define the {\em open ball} $B_\varepsilon(x)$ by
\[B_\varepsilon(x)=\{y\in X : d(x,y) <\varepsilon\},\]
and say that $U\subseteq X$ is {\em open} if, for all $x\in U$, there is $\varepsilon > 0$ such that $B_\varepsilon (x) \subseteq U$. The Euclidean spaces $\mathbb R^n$ are all examples of metric spaces and are always assumed to be equipped with the topology we have just defined, as is the set of rational numbers, denoted $\mathbb Q$. The set of open balls on a metric space is a basis for its topology, in the following sense:

\begin{defn}\label{DefBasis}
A collection $\mathcal{B}$ of subsets of a set $X$ is a {\em basis} if
\begin{enumerate}
\item\label{ItBasisOne} $\bigcup_{B\in\mathcal{B}}B=X$, and
\item\label{ItBasisTwo} whenever $B_0,B_1\in\mathcal{B}$ and $x\in B_0\cap B_1$, there exists $B_2\subseteq B_0\cap B_1$ such that $x\in B_2.$
\end{enumerate}
The basis $\mathcal B$ {\em generates} a topology on $X$ defined by letting $U\subseteq X$ be open if and only if, for every $x
\in U$, there is $B\in \mathcal B$ with $x\in B \subseteq U$.
\end{defn}

Topological spaces can also be seen as a generalization of preorders. If $\mathfrak W$ is a preordered set, consider the topology $\mathcal D_\peq$ on $|\mathfrak W|$ given by setting $U\subseteq|\mathfrak W|$ to be open if and only if, whenever $w\in U$, we have ${\downarrow} w\subseteq U$ (so that the sets of the form ${\downarrow} w$ provide a basis for $\mathcal D_\peq$). We call $\mathcal D_\peq$ the {\em down-set topology of $\peq$.} Similarly, the family of upwards-closed subsets of $|\mathfrak W|$ will be denoted by $\mathcal U_\peq$, and is the {\em up-set topology of $\peq$.} Topologies of this form are Alexandroff topologies \cite{aleksandroff}:

\begin{defn}\label{DefAlex}
A topological space $\mathfrak A$ is an {\em Alexandroff space} if any one of the following equivalent conditions occurs:
\begin{enumerate}

\item whenever $\mathcal O\subset \mathcal T_\mathfrak A$, then $\bigcap \mathcal O\in \cl T_\fr A$;

\item every $x\in |\fr A|$ has a $\subseteq$-least neighborhood;

\item\label{ItAlexTwo} there is a preorder $\peq$ on $|\mathfrak A|$ such that $\cl T_\fr A=\cl D_{\peq}$, or

\item there is a preorder $\peq$ on $|\mathfrak A|$ such that $\cl T_\fr A=\cl U_{\peq}$.

\end{enumerate}
\end{defn}

It is readily verified that if $\cl T_\fr A=\cl D_{\peq}$, then $\peq$ is uniquely defined,\footnote{In this article we will be biased towards using the down-set topology since it is more natural for the constructions we will introduce in \S \ref{SecFour}.} and we will denote it by $\peq_\fr A$. We remark that intuitionistic logic cannot distinguish between preordered sets and partially ordered sets, so without loss of generality we may work only with Alexandroff spaces generated by a poset. We will tacitly identify $(|\fr A|,\peq_\fr A)$ with $(|\fr A|,\cl T_\fr A)$.

Observe that all finite topological spaces are Alexandroff; in fact, all {locally} finite spaces are Alexandroff. The following is easily verified:

\begin{lem}
Say that a topological space $\fr A$ is {\em locally finite} if every point of $|\fr A|$ has a neighborhood $U$ such that $\# U < \omega$. Then, every locally finite space is Alexandroff.
\end{lem}

It will also be useful to characterize the continuous functions on a poset:
\begin{lem}
If $\mathfrak W,\fr V$ are preorders and $g\colon |\mathfrak W|\to |\mathfrak V|$, then $g$ is continuous with respect to the down-set topologies on $\mathfrak W,\fr V$ if and only if, whenever $v\peq_\fr W w$, it follows that $g(v)\peq_\fr V g(w)$.
\end{lem}

In other words, continuous maps on posets are simply monotone maps. Alexandroff spaces will be very useful to us throughout the text, as our goal is to reduce intuitionistic temporal logic over arbitrary dynamical systems to a computably bounded set of finite structures. With this, we are ready to define the logics that we are interested in.

\section{Syntax and Semantics}\label{SecBasic}

In this section we will introduce the logic $\icltl$ and its semantics; in fact, we will first define an extended logic $\itlcplus$, containing both our logic and Kremer's logic $\itlk$.
Fix a countably infinite set $\mathbb P$ of propositional variables. Then, the {\em full (intuitionistic temporal) language} $\lanfull$ is defined by the grammar (in Backus-Naur form)
\[\varphi,\psi := \ \  \bot \  | \   p  \ |  \ \varphi\wedge\psi \  |  \ \varphi\vee\psi  \ |  \ \varphi\imp \psi  \ |  \ \nx\varphi \  | \  \ps\varphi \  |  \ \nec\varphi  \ |  \ \exists\varphi \  | \  \forall \varphi, \]
where $p\in \mathbb P$. Here, $\nx$ is read as `next', $\ps$ as `eventually', and $\nec$ as `henceforth'. We also use $\ineg\varphi$ as a shorthand for $\varphi\imp \bot$ and $\varphi\iiff \psi$ as a shorthand for $(\varphi\imp \psi) \wedge (\psi\imp\varphi)$.

Denote the set of subformulas of $\varphi\in\lanfull$ by ${\mathrm{sub}}(\varphi)$, and the length of $\varphi$ (or, more precisely, $\# {\mathrm{sub}}(\varphi)$) by $\lgt\varphi$. We will also be interested in certain sublanguages of $\lanfull$, which only allow some set of modalities $M\subseteq \{\nx,\ps,\nec,\exists,\forall\}$, and we will denote such a fragment by $\lanfull\upharpoonright M$. Note that $\lanfull \upharpoonright M$ always contains Booleans and implication, even though they will not be listed in $M$. We will be primarily interested in the language $\landi=\lanfull\upharpoonright \{\circ,\ps,\forall\}$.
Formulas of $\lanfull$ are interpreted on dynamical systems over topological spaces, or {\em dynamical topological systems.}

\begin{defn}\label{DefSem}
A {\em dynamical (topological) system} is a triple $\mathfrak X=(|\mathfrak X|,\mathcal{T}_\mathfrak X,f_\mathfrak X)$
where $(|\mathfrak X|,\mathcal{T}_\mathfrak X)$ is a topological space and $f_\fr X:|\mathfrak X|\to |\mathfrak X|$ is a continuous function.
A {\em valuation on $\mathfrak X$} is a function $\lb\cdot\rb\colon\lanfull \to \mathcal T_\mathfrak X$ such that
\[
\begin{array}{rcl}
\lb\bot\rb&=&\varnothing \\
\lb\varphi\wedge\psi\rb &=&\lb\varphi\rb\cap \lb\psi\rb\\
\lb\varphi\vee\psi\rb &=&\lb\varphi\rb\cup \lb\psi\rb\\
\lb\varphi\imp\psi\rb &=&\big ( (|\mathfrak X|\setminus\lb\varphi\rb)\cup \lb\psi\rb\big )^\circ\\
\val{\nx\varphi}&=&f^{-1}_\mathfrak X\val\varphi\\
\val{\ps\varphi}&=& \bigcup_{n<\omega}f^{-n}_\mathfrak X\val\varphi\\
\val{\nec\varphi}&=&\big (\bigcap_{n<\omega}f^{-n}\val\varphi\big)^\circ\\
\val{\exists\varphi}&=&
|\mathfrak X|\text{ if $\val\varphi\not=\varnothing$ and
$\varnothing$ otherwise;}\\
\val{\forall\varphi}&=&
|\mathfrak X|\text{ if $\val\varphi=|\mathfrak X|$ and
$\varnothing$ otherwise.}
\end{array}
\]
A dynamical system $\mathfrak X$ equipped with a valuation $\lb\cdot\rb_\mathfrak X$ is a {\em (dynamical topological) model}.
\end{defn}

If $f_\mathfrak X$ is a homeomorphism we will say that $\mathfrak X$ is {\em invertible.}
Note that, for a propositional variable $p$, $\lb p\rb$ can be any subset of $|\mathfrak X|$, {\em provided} it is open. The rest of the clauses are standard from either intuitionistic or temporal logic, with the exception of $\nec\varphi$. The interpretation of the latter is due to Kremer, and we will discuss it further in the next section.

As a general convention, we let structures inherit the properties of their components, so that for example an Alexandroff model is a model $\fr A$ such that $(|\fr A |,\cl T_\fr A)$ is an Alexandroff space. In the setting of Alexandroff models, the semantics for implication simplifies somewhat.

\begin{lem}\label{LemmAlex}
Let $\mathfrak A$ be a model based on an Alexandroff space, and $\varphi$ be any formula. Then,
\begin{enumerate}

\item

If $\cl T_\fr A=\cl D_{\peq}$, $x\in |\fr A|$ and $\varphi,\psi\in \lanfull$, then $x\in \val{\varphi\imp \psi}_\fr A$ if and only if, for all $y\peq  x$, if $y\in \val {\varphi }_\fr A$, then $y\in \val {\psi }_\fr A$.

\item

Similarly, if $\cl T_\fr A=\cl U_{\peq}$, $x\in |\fr A|$ and $\varphi,\psi\in \lanfull$, then $x\in \val{\varphi\imp \psi}_\fr A$ if and only if, for all $y\cca x$, if $y\in \val {\varphi }_\fr A$, then $y\in \val {\psi }_\fr A$.

\item $\val{\nec\varphi}_\fr A=\bigcap_{n<\omega } f^{-n}\val\varphi_\fr A.$

\end{enumerate}

\end{lem}

\proof
For the first claim, observe that if $x\in E\subseteq |\fr A|$, then $x\in E^\circ$ if and only if $\mathop\downarrow x\subseteq E$. From this and the definition of $\val{\varphi \imp \psi}_\fr A$, the claim follows. The second claim is analogous, but replacing $\peq$ by $\cca$. For the third, since infinite intersections of open sets are open in an Alexandroff space, it follows that $\bigcap_{n<\omega } f^{-n}\val\varphi_\fr A$ is open, and thus
\[\val{\nec\varphi}_\fr A=\big(\bigcap_{n<\omega } f^{-n}\val\varphi_\fr
  A\big)^\circ=\bigcap_{n<\omega } f^{-n}\val\varphi_\fr A.
  \tag*{\qedhere}
\]
\endproof

Validity is then defined in the usual way:

\begin{defn}
Given a model $\mathfrak X$ and a formula $\varphi\in \lanfull$, we say that $\varphi$ is {\em valid} on $\mathfrak X$, written $\mathfrak X\models\varphi$, if $\val\varphi_\mathfrak X =|\mathfrak X|$. If $\mathfrak X$ is a dynamical system, we write $\mathfrak X\models\varphi$ if $(\mathfrak X,\val\cdot)\models \varphi$ for every valuation $\val\cdot$ on $\mathfrak X$. If $\Omega$ is a class of dynamical systems or models, we say that $\varphi\in\lanfull$ is {\em valid on $\Omega$} if, for every $\mathfrak X\in \Omega$, $\mathfrak X\models\varphi$. If $\varphi$ is not valid on $\Omega$, it is {\em falsifiable on $\Omega$.} 

We define the logics $\icltl$ and $\itlcplus$ to be the set of formulas of $\landi$ and $\lanfull$, respectively, that are valid over the class of all dynamical systems.
\end{defn}

One of the key differences between classical and intuitionistic logic has to do with the fact that any formula is classically equivalent to its double negation, but intuitionistically this need not be the case. To this end, it is convenient to elucidate the topological meaning of double negation.

\begin{lem}\label{LemmDoubleNeg}
Let $\mathfrak X$ be any model and $\varphi\in \lanfull$. Then, $\val{\ineg\ineg\varphi} = \big ( \overline{\val\varphi} \big )^\circ$; that is, $\ineg\ineg\varphi$ is true precisely on the interior of the closure of $\val\varphi$.
\end{lem}

We leave the proof to the reader; it can be a good exercise to familiarize oneself with the topological semantics of intuitionistic logic. As a corollary, we obtain the following fact, which will be useful throughout the text: given formulas $\varphi,\psi$ and a model $\mathfrak X$, we have that $\mathfrak X\models \varphi\imp \ineg\ineg \psi$ if and only if $\val\psi$ is {\em dense} in $\val \varphi$; that is, if $\val \varphi \subseteq \overline{\val \psi}$.

\begin{exam}
Consider a model $(\mathbb R,f,\val\cdot)$, where $\mathbb R$ is the real line equipped with the standard topology, and $\val p = (-1,0) \cup (0,1)$; the value of $f$ is unimportant for this example, and we may take e.g.~$f(x) = x$. Then, in view of Lemma \ref{LemmDoubleNeg},
\[\val{\ineg\ineg p} = \Big ( \overline{(-1,0) \cup (0,1)}\Big)^\circ =  [-1,1]  ^\circ = (-1,1).\]
In particular, $0\not\in \val {\ineg \ineg p\imp  p}$, showing that the latter formula is not intuitionistically valid.
\end{exam}

Next we exhibit a principle valid in Alexandroff systems that, as we will see later, is not valid in general.

\begin{exam}\label{ExAlekx}
The formula 
\[\forall(\ineg p \vee \ps p) \Rightarrow \ineg \ps p \vee \ps p\]
is valid on any model based on an Alexandroff space. To see this, suppose that $\mathfrak A$ is any such model with $\mathcal T_\mathfrak A=\mathcal D_\peq$, and suppose that $w\in \val{\forall(\ineg p\vee \ps p)}_\mathfrak A$, so that $\val{\ineg p \vee \ps p}_\fr A = |\mathfrak A|$. If $w \in \val{\ineg\ps p}_\mathfrak A$ we are done, otherwise there is $v\acc w$ such that $v\in \val{\ps p}_\mathfrak A$, from which we obtain $f^n_\mathfrak A(v)\in \val{p}_\mathfrak A$ for some $n$.

Note that $f^n_\mathfrak A(v) \acc f^n_\mathfrak A(w)$, from which it follows that  $f^n_\mathfrak A(w) \not \in \val{\ineg p}_\mathfrak A$. But $f^n_\mathfrak A(w)  \in \val{\ineg p \vee \ps p}_\mathfrak A$, so that $f^n_\mathfrak A(w)\in \val {\ps p}_\mathfrak A$ and thus $f^{n+m}_\mathfrak A(w)\in \val {p}_\mathfrak A$ for some $m\geq 0$, which implies that $w \in \val {\ps p}_\mathfrak A$, as required.
\end{exam}
However, as we will see later, the above formula $\varphi$ is not valid over the class of {\em all} dynamical systems. Next, let us show that our full language can be simplified somewhat; in particular, it admits $\exists$-elimination.

\begin{lem}\label{LemmExDefin}
Given any model $\mathfrak X$ and any formula $\varphi$, $\val{\exists\varphi}_\mathfrak X = \val{\ineg \forall \ineg \varphi}_\mathfrak X.$
\end{lem}

\proof
Note that $\val{\exists\varphi}_\mathfrak X \in \{\varnothing,|\mathfrak X|\}$. If $\val{\exists\varphi}_\mathfrak X = \varnothing$, then $\val\varphi_\mathfrak X=\varnothing$, and thus $\val{\ineg\varphi}_\mathfrak X=|\mathfrak X|$. It follows that $\val{\forall\ineg\varphi}_\mathfrak X=|\mathfrak X|$, and thus $\val{\ineg\forall\ineg\varphi}_\mathfrak X=\varnothing$.

Otherwise, $\val{\exists\varphi}_\mathfrak X = |\mathfrak X|$, which means that $\val\varphi_\mathfrak X\not=\varnothing.$ Hence $\val{\ineg\varphi}_\mathfrak X\not = |\mathfrak X|$; but then, $\val{\forall  \ineg \varphi}_\mathfrak X = \varnothing$, therefore $\val{\ineg \forall  \ineg \varphi}_\mathfrak X = |\mathfrak X|$.
\endproof

Thus we may turn our attention to $\lanfull\upharpoonright \{\nx,\ps,\nec,\forall\}$, or, generally speaking, to languages without $\exists$; however, we do want $\exists$ to be definable, as it will allow us to capture minimality (see \S \ref{SecMinimal}). In contrast, if $\lanfull_\exists$ denotes the $\forall$-free fragment $\lanfull\upharpoonright \{\nx,\ps,\nec,\exists \}$, then $\forall$ is not $\lanfull_\exists$-definable, as we will see next.
Below, if $X$ is any set, $f\colon X\to X$ and $U\subseteq X$, we say that $U $ is {\em $f$-invariant} if $f(U)\subseteq U$, and we say that two valuations $\val\cdot$, $\val\cdot'$ {\em agree} on $U$ if $\val p \cap U = \val p '$ for all propositional variables $p$.

\begin{lem}\label{lemExistsDense}
Let $\mathfrak X$ be a dynamical system and let $\val\cdot$, ${\val\cdot}'$ be two valuations on $\mathfrak X$ that agree on a dense, open, $f_\mathfrak X$-invariant set $U \subseteq |\mathfrak X|$.
Then $\val \varphi \cap U = \val \varphi' \cap U$ for all formulas $\varphi \in \lanfull_\exists$.
\end{lem}

\proof
The proof proceeds by a standard structural induction on $\varphi$. We treat only the cases where $\varphi$ is of the form $\nec \psi$ or $\exists \psi$ as examples.

If $x\in \val{\nec \psi}$, there is a neighbourhood $V$ of $x$ such that for all $y\in V$ and $n<\omega$, $f^n_\mathfrak X(y) \in \val \psi$. Since $U$ is open we may assume that $V\subseteq U$, in which case if $y\in V$ and $n<\omega$ we obtain $f^n_\mathfrak X(y) \in U$ from the assumption that $U$ is $f_\mathfrak X$-invariant. Thus we may apply the induction hypothesis to see that $f^n_\mathfrak X(y) \in \val \psi '$, so that $V$ witnesses that $x\in \val {\nec \psi}'$, as needed. The converse implication is entirely symmetrical.

Finally, if $\varphi = \exists \psi$ and $x\in \val{\exists \psi}$, then $\val\psi \not= \varnothing$. Since $U$ is dense there is some $y\in \val \psi \cap U$, which by induction hypothesis gives us $y\in \val \psi '$ and thus $x\in \val {\exists \psi}'$.
\endproof

From this it is straightforward to construct two models that agree on $\lanfull_\exists$ but not on $\lanfull$.

\begin{prop}
There is no formula in $\lanfull_\exists$ equivalent to $\forall p$ over the class of finite, invertible dynamical systems.
\end{prop}

\proof
Let $\mathfrak X = (W,\mathcal D_{\leq}, {\bm 1}_W)$ where $W=\{0,1\}$, $ 0 \leq 1$ and ${\bm 1}_W$ is the identity function on $W$. Let $\val\cdot$, $\val \cdot'$ be valuations so that $\val p= W$, $\val p ' = \{0\}$ and $\val q = \val q ' = \varnothing$ for all $q\not=p$. Then $ \val \cdot$ and $\val \cdot '$ agree on $\{0\}$, which is dense, open, and trivially ${\bm 1}_W$-invariant. It follows that for all $\varphi \in \lanfull_\exists$, $0 \in \val \varphi$ if and only if $0 \in \val \varphi '$; however $0 \in \val {\forall p}$ but $0\not \in \val {\forall p}'$, so that $\varphi$ is not equivalent to $\forall p$.
\endproof

Similarly, $\ps$ and $\nec$ cannot be defined in terms of each other, even over the class of finite dynamical systems \cite{Balbiani2017}, but we omit $\nec$ from our language for technical reasons that will be discussed later. Meanwhile, we have included $\forall$ in $\landi$ but not $\exists$; since the latter is definable, we lose nothing by omitting it.

\section{Related Logics}\label{SecComp}

Before we continue our analysis of $\icltl$, let us mention some related systems. First we observe that $\itlcplus$ is an extension of the logic $\itlk$ introduced by Kremer in \cite{KremerIntuitionistic}. To be precise, Kremer considers the language\footnote{In order to match the presentation in \cite{KremerIntuitionistic}, we do not include $\forall$ in $\lank$. However, it is an inessential choice for our purposes, and adding $\forall$ to $\lank$ would not affect our discussion.} $\lank=\lanfull\upharpoonright \{\nx,\nec\}$. There he also presents the next example, which will help elucidate the semantics of $\nec$.
\begin{exam}\label{ExKremer}
Consider the dynamical system $(\mathbb R, f)$ where $f \colon \mathbb R \to \mathbb R$ is given by
\[
f(x)=
\begin{cases}
0 & \text{if $x\leq 0$,}\\

2x &\text{otherwise.}
\end{cases}
\]
Suppose that $\val p= (-\infty,1)$.

Observe that, for all $n$, $f^n(0)=0\in \val p$. Thus $0\in \bigcap_{n<\omega}f^{-n}\val p$. Moreover, if $x<0$, then $x\in \val p$ and $f(x)=0$, so once again $x\in \bigcap_{n<\omega}f^{-n}\val p$.

On the other hand, if $x>0$, then  $f^n(x)=2^n x$. Since $2^n \to\infty$ as $n\to \infty$ it follows that $2^n x> 1$ for $n$ large enough, and thus $f^n(x)\not\in \val p$. It follows that $x\not\in \bigcap_{n<\omega}f^{-n}\val p$.
But then
\[\bigcap_{n<\omega}f^{-n}\val p =(-\infty,0],\]
which is not open. Since intuitionistic truth values must always be open, we need to ``approximate'' it by an open set; this is done in the standard way by taking the interior, and thus
\[\val{\nec p}=\big (\bigcap_{n<\omega}f^{-n}\val p \big)^\circ =(-\infty,0).\]
\end{exam}

Kremer used Example \ref{ExKremer} in \cite{KremerIntuitionistic} to prove the that several key principles of $\sf LTL$ fail intuitionistically.

\begin{thm}\label{TheoKrem}
Let ${\sf ITL}^\mathbb R$ be the set of formulas of $\lanfull$ valid over the class of dynamical systems based on $\mathbb R$ (with the usual topology). Then, $\nec p\imp \nx\nec p$, $\nec\nx p\imp \nx\nec p$, and $\nec p\imp \nec\nec p \not \in {\sf ITL}^\mathbb R$.
\end{thm}

This shows that some care must be taken when working with $\nec$, and as we will discuss in Example \ref{ExKremerTwo}, the decidability proof we will present would not go through if we included it (at least not without some non-trivial modifications). Note, however, that in view of Lemma \ref{LemmAlex}, this issue disappears when working over the class of Alexandroff spaces.

Moreover, Alexandroff spaces will allow us to compare our semantics with the expanding frames presented in Boudou et al.~\cite{Boudou2017}. There, we use relational structures $(W,{\peq},f)$, where $w\peq v$ implies that $f(w)\peq f(v)$, and let $w\in W$ satisfy $\varphi\imp \psi$ if, whenever $v\cca w$, if $v$ satisfies $\varphi$, it also satisfies $\psi$; note that these are precisely the truth conditions with respect to the up-set topology $\mathcal U_\peq$. Structures with these properties are similar to expanding products in modal logic \cite{pml}, and we will denote the logic of such frames over $\landi$ and $\lanfull$ by $\itle$ and $\itlepl$, respectively. We proved in \cite{Boudou2017} that $\itlepl$ is decidable, with a similar super-exponential bound as we will obtain for $\icltl$, although the techniques are quite different.\footnote{Note that in \cite{Boudou2017} we do not include a universal modality, but the proof of decidability can be readily modified to accommodate it.} This is proven by showing that $\itlepl$ has the effective finite model property; as we will see, the same does not hold for $\icltl$, which is why we need to work instead with quasimodels (see \S \ref{SecNDQ}).

In \cite{Boudou2017}, we also consider frames with the additional `backward confluence' property ${\acc}\circ f\subseteq f\circ{\acc}$; these systems are also considered in \cite{KojimaNext}, there called {\em functional Kripke frames.} We will follow \cite{Boudou2017} and call them {\em persistent frames,} since they are related to standard (persistent) products of modal logics \cite{mdml}. The intuitionistic temporal logic of persistent frames is denoted $\itlaipl$.

Backwards confluence may also be written as $\mathord \uparrow f(w)\subseteq f[\mathord \uparrow w]$; in other words, $f$ is an open map.
Hence, expanding frames are precisely the Alexandroff dynamical systems, and persistent frames are the Alexandroff dynamical systems with an open map. It is immediate that $\itlcplus\subseteq{\itlepl}\subseteq \itlaipl$; in fact, both inclusions are strict. In Balbiani et al.~\cite{Balbiani2017} it is shown that $(\circ p\imp \circ q) \imp \circ (p\imp q) \in \itlaipl\setminus \itlepl$, and later we will show that $\icltl \ne {\itle }$.

Finally, it will be convenient to compare $\itlcplus$ to dynamic topological logic. Since the latter's base logic is classical, we may use a simpler syntax for $\sf DTL$, using the language $\lanclass$ given by the grammar
\[ \bot \  | \   p  \ |  \ \varphi\to\psi  \ |  \ \blacksquare \varphi \  | \  \nx\varphi \  |  \ \nec\varphi \  |  \ \forall \varphi. \]
We can then define $\neg,\wedge,\vee,\blacklozenge,\ps$ using standard classical validities. Since we will go back and forth between the classical and intuitionistic interpertations, we will always use $\neg,\to$ for classical logics and $\ineg,\imp$ for intuitionistic ones to avoid confusion.

\begin{table}

\begin{center}

\begin{tabular}{|c|c|c|}
\hline
Language & Connectives & Modalities \\
\hline
\hline
$ \phantom{\displaystyle{a \choose b}} \lanfull  \phantom{\displaystyle{a \choose b}}$ & $\bot,\vee,\wedge,\imp$ & $\nx, \ps, \nec,\exists,\forall$\\
\hline
$ \phantom{\displaystyle{a \choose b}} \landi  \phantom{\displaystyle{a \choose b}}$ & $\bot,\vee,\wedge,\imp$ & $\nx, \ps, \forall$\\
\hline
$ \phantom{\displaystyle{a \choose b}} \lank  \phantom{\displaystyle{a \choose b}}$ & $\bot,\vee,\wedge,\imp$ & $\nx, \nec$\\
\hline
$ \phantom{\displaystyle{a \choose b}} \lanfull_\exists  \phantom{\displaystyle{a \choose b}}$ & $\bot,\vee,\wedge,\imp$ & $\nx, \ps, \nec,\exists$\\
\hline
$ \phantom{\displaystyle{a \choose b}} \lanclass  \phantom{\displaystyle{a \choose b}}$ & $\bot,\rightarrow$ & $\nb, \nx, \nec, \forall$\\
\hline
\end{tabular}

\end{center}

\caption{Formal languages considered in the text. Note that, in order to simplify notation, we indicate at most one distinctive modality.}

\end{table}

Given a dynamical system $\mathfrak X$, a {\em classical valuation on $\mathfrak X$} is a function $\cval\cdot\colon\lanclass\to \mathcal P(|\mathfrak X|)$ such that
\[
\begin{array}{lcllcl}
\cval \bot&=&\varnothing &
\cval{\varphi\to\psi}&=& (|\mathfrak X|\setminus\cval\varphi)\cup \cval\psi\\[.8ex]
\cval{\blacksquare\varphi}&=&(\cval\varphi)^\circ&
\cval{\nx\varphi}&=&f^{-1}_\mathfrak X\cval\varphi\\[.8ex]
\cval{\nec\varphi}&=&\displaystyle\bigcap_{n<\omega}f^{-n}_\mathfrak X\cval\varphi\phantom{aaaaa}&
\cval{\forall\varphi}&=&
|\mathfrak X|\text{ if $\cval\varphi=|\mathfrak X|$ and
$\varnothing$ otherwise.}\\
\end{array}
\]
Note that valuations of formulas are no longer restricted to open sets. Our intuitionistic temporal logic may then be interpreted in $\sf DTL$ via the G\"odel-Tarski translation $\cdot^\blacksquare$, defined as follows:

\begin{defn}
Given $\varphi\in \lanfull$, we define $\varphi^\blacksquare\in \lanclass$ recursively by setting
\[
\begin{array}{lcllcl}
\bot^\nb & = & \bot &
p^\nb & = & \nb p \\
(\varphi \odot  \psi )^\nb & = & \varphi^\nb \odot \psi^ \nb \phantom{aaaaa}&
(\varphi\imp \psi )^\nb & = &\nb ( \varphi^\nb \to \psi^ \nb ) \\
(\nec \varphi)^\nb & = &\nb \nec \varphi^\nb  &
(\boxdot  \varphi)^\nb & = & \boxdot \varphi^\nb ,
\end{array}
\]
where $\odot\in \{\wedge,\vee\}$ and $\boxdot \in \{\nx,\ps,\exists,\forall\}$.
\end{defn}

For example,
\[( \ps  \nec p \imp \bot  )^\blacksquare= \blacksquare (\ps \nb  \nec \nb p \to \bot  ).\]
The following can then be verified by a simple induction on $\varphi$:

\begin{lem}\label{LemGT}
Let $\varphi\in \lanfull$, and $\mathfrak X$ be any dynamic topological system. Suppose that $\val\cdot$ is an intuitionistic valuation and $\cval\cdot$ a classical valuation such that, for every atom $p$, $(\cval{p})^\circ=\val p$. Then, for every formula $\varphi$, $\val\varphi=\cval{\varphi^\blacksquare}$.
\end{lem}

From Lemma \ref{LemGT} and the fact that $\sf DTL$ is computably enumerable \cite{me2}, we immediately obtain the following.

\begin{thm}\label{TheoStarCE}
The logic $\itlcplus$ is computably enumerable.
\end{thm}

Note that \cite{me2} considers a language without a universal modality, but adding it does not affect computable enumerability. Moreover, we proved in \cite{me:metric} that any satisfiable $\lanclass$ formula was also satisfiable on a model based on $\mathbb Q$; as before, the universal modality can readily be incorporated, giving us the following:

\begin{thm}\label{TheoQComplete}
Let ${\sf ITL}^{\mathbb Q}$ be the set of $\lanfull$-formulas that are valid over the class of dynamical systems based on the rational numbers with the usual topology. Then, ${\sf ITL}^{\mathbb Q}=\itlcplus$.
\end{thm}

However, $\sf DTL$ is undecidable \cite{konev}, and hence Theorem \ref{TheoStarCE} does not settle the decidability of $\itlcplus$. It is not known whether the full $\itlcplus$ is decidable; however, in the rest of this article we will show that this is indeed the case for $\icltl$.

\section{Labelled Systems}\label{SecNDQ}

Our decidability proof is based on (non-deterministic) quasimodels, introduced in \cite{me2}. In this section we will introduce labeled systems, which generalize both quasimodels and dynamical topological models, thus allowing us to view both of them in a unified framework. Given a set $\Sigma\subset\landi$ that is closed under subformulas, we say that $\Phi\subseteq \Sigma$ is a {\em $\Sigma$-type} if the following occur.
\begin{enumerate}

\item $\bot\not\in \Sigma$.

\item If $\varphi\wedge\psi\in \Sigma$, then $\varphi\wedge\psi\in \Phi$ if and only if $\varphi,\psi\in \Phi$.

\item If $\varphi\vee\psi\in \Sigma$, then $\varphi\vee\psi\in \Phi$ if and only if $\varphi\in\Psi$ or $\psi\in \Phi$.

\item If $\varphi\imp \psi\in \Sigma$, then

\begin{enumerate}

\item  $\varphi\imp \psi\in\Phi$ implies that $\varphi\not \in \Phi$ or $\psi\in\Phi$ and

\item $\psi \in \Phi$ implies that $\varphi\imp \psi \in\Phi$.

\end{enumerate}

\item If $\ps\psi\in \Sigma$ and $\psi\in \Phi$, then $\ps\psi\in \Phi$.

\end{enumerate}
The set of $\Sigma$-types will be denoted by $\type\Sigma$.
Many times we want $\Sigma$ to be finite, in which case given $\Delta\subseteq \landi$ we write $\Sigma\Subset \Delta$ to indicate that $\Sigma\subseteq \Delta$ is finite and closed under subformulas.
Note that $\type\Sigma$ is partially ordered by $\subseteq$, and we will endow it with the up-set topology $\mathcal U_\subseteq$. For $\Phi\in\type\Sigma$, say that a formula $\varphi\imp \psi\in\Sigma$ is a {\em defect} of $\Phi$ if $\varphi \imp \psi\not\in\Phi$ and $\varphi\not\in \Phi$. The set of defects of $\Phi$ will be denoted $\partial_\Sigma\Phi$, or simply $\defect\Phi$ when $\Sigma$ is clear from context.

\begin{defn}\label{frame}
Let $\Sigma\subset\landi$ be closed under subformulas.
We say that a {\em $\Sigma$-labeled space} is a triple $\mathfrak W= ( |\mathfrak W|,\mathcal T_\mathfrak W,\ell_\mathfrak W )$, where $( |\mathfrak W|,\mathcal T_\mathfrak W )$ is a topological space and $\ell_\mathfrak W\colon | \mathfrak W | \to \type\Sigma$ a continuous function such that for all $w\in |\mathfrak W|$, whenever $\varphi\imp \psi\in \defect{ \ell_\mathfrak W(w)}$ and $U$ is any neighborhood of $w$, there is $v\in U$ such that $\varphi\in \ell_\mathfrak W(v)$ and $\psi\not\in \ell_\mathfrak W(v)$. Such a $v$ {\em revokes} $\varphi\imp \psi$.

The $\Sigma$-labeled space $\mathfrak W$ {\em falsifies} $\varphi\in\mathcal L$ if $\varphi\in \Sigma\setminus \ell_\mathfrak W(w)$ for some $w\in |\mathfrak W|$, and $\ell_\mathfrak W$ is {\em honest} if, for every $w\in |\mathfrak W|$ and every $\forall\varphi\in \Sigma$, we have that $\forall\varphi\in \ell_\mathfrak W(w)$ if and only if $\varphi\in \ell_\mathfrak W(v)$ for every $v\in |\mathfrak W|$.
\end{defn}

If $\mathfrak W$ is a labelled space, elements of $|\mathfrak W|$ will sometimes be called {\em worlds.} As usual, we may write $\ell$ instead of $\ell_\mathfrak W$ when this does not lead to confusion. Since we have endowed $\type\Sigma$ with the topology $\mathcal U_\subseteq$, the continuity of $\ell$ means that for every $w\in|\mathfrak W|$, there is a neighborhood $U$ of $w$ such that, whenever $v \in U$, $\ell_\mathfrak W(w)\subseteq \ell_\mathfrak W(v)$.

Note that not every subset $U$ of $|\mathfrak W|$ gives rise to a substructure that is also a labelled space; however, this is the case when $U$ is open.

\begin{lem}\label{LemmOpenSubst}
If $\Sigma\subseteq \landi$ is closed under subformulas, $\mathfrak W$ is a $\Sigma$-labelled space, and $U\subseteq |\mathfrak W|$ is open, then $\mathfrak W\upharpoonright U$ is a $\Sigma$-labelled space.
\end{lem}

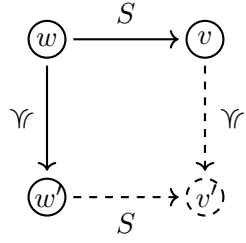
\begin{figure}

\begin{center}

\begin{tikzpicture}[scale=0.7]

\draw[thick] (0,0) circle (.35);

\draw (0,0) node {$w$};

\draw[thick,->] (.5,0) -- (2.5,0);

\draw (1.5,.5) node {$S$};

\draw[thick] (3,0) circle (.35);

\draw (3,.02) node {$v$};

\draw[thick] (0,-3) circle (.35);

\draw (.04,-3+.03) node {$w'$};

\draw[thick,->,dashed] (.5,-3) -- (2.5,-3);

\draw (1.5,-3.5) node {$S$};

\draw[thick,<-] (0,-2.5) -- (0,-.5);

\draw (-.5,-1.5) node {{\large$\rotatebox[origin=c]{90}{$\peq$}$}};

\draw[thick,dashed] (3,-3) circle (.35);

\draw (3+.04,-3+.03) node {$v'$};

\draw[thick,<-,dashed] (3,-2.5) -- (3,-.5);

\draw (3.5,-1.5) node {{\large$\rotatebox[origin=c]{90}{$\peq$}$}};

\end{tikzpicture}

\end{center}

\caption{If $(A,{\peq})$ is a partial order and $S \subseteq A \times A$, then $S$ is continuous with respect to the down-set topology $\mathcal D_\peq$ whenever the above diagram can always be completed.}
\label{FigContRel}

\end{figure}

For our purposes, a {\em continuous relation} on a topological space is a relation under which the preimage of any open set is open (note that this is not the standard definition of a contiuous relation, which is more involved). In the context of an Alexandroff space with the down-set topology, a continuous relation $S$ is one that satisfies the forward confluence property depicted in Figure \ref{FigContRel}: if $w'\peq w \mathrel S v$, then there is $v'$ such that $w'\mathrel S v' \peq v$.

\begin{defn}\label{compatible}
Let $\Sigma\subset\landi$ be closed under subformulas. Suppose that $\Phi,\Psi\in\type\Sigma$. The ordered pair $(\Phi,\Psi)$ is {\em sensible} if
\begin{enumerate}
\item for all $\nx\varphi\in \Sigma$, $\nx\varphi\in \Phi$ if and only if $ \varphi\in \Psi$,
\item for all $\ps\varphi\in \Sigma$, $\ps\varphi\in \Phi$ if and only if $\varphi\in\Phi$ or $\ps\varphi\in \Psi$, and
\item for all $\forall \varphi\in \Sigma$, $\forall \varphi\in \Phi$ if and only if $\forall\varphi\in \Psi$.
\end{enumerate}
Likewise, a pair $(w,v)$ of worlds in a labelled space $\mathfrak W$ is sensible if $(\ell (w),\ell (v))$ is sensible.
A continuous relation
$S\subseteq |\mathfrak W|\times |\mathfrak W|$
is {\em sensible} if every pair in $S$ is sensible.
Further, $S$ is {\em $\omega$-sensible} if it is serial and, whenever $\ps\varphi\in \ell(w)$, there are $n\geq 0$ and $v$ such that $w \mathrel S^n v$ and $\varphi\in \ell(v)$.

A {\em labelled system} is a labelled space $\mathfrak W$ equipped with a sensible relation $S_\mathfrak W\subseteq |\mathfrak W|\times|\mathfrak W|$; if moreover $\ell_\mathfrak W$ is honest and $S_\mathfrak W$ is $\omega$-sensible, we say that $\mathfrak W$ is a {\em well $\Sigma$-labelled system.}

\end{defn}

Any dynamic topological model can be regarded as a well-labelled system. If $\Sigma\subseteq \landi$ is closed under subformulas, $\mathfrak{X}$ is a model and $x\in |\mathfrak X|$, we can assign a $\Sigma$-type $\ell_\mathfrak X(x)$ to $x$ given by
\[\ell_\mathfrak X(x)=\cbra\psi\in \Sigma :x\in \val\psi_\mathfrak X \cket.\]
We also set $S_\mathfrak X=f_\mathfrak X$; it is obvious that $\ell_\mathfrak X$ is honest and $S_\mathfrak X$ is $\omega$-sensible.
Henceforth we will tacitly identify $\mathfrak X$ by its associated well $\landi$-labelled system.

\begin{exam}\label{ExKremerTwo}
Example \ref{ExKremer} can be used to show that there are some difficulties when treating $\nec$ using labelled systems. In that example $\val{\nec p}=(-\infty,0)$, so that, for instance, $-1\in \val{\nec p}$, but $f(-1)=0 \not \in \val{\nec p}$. If we wanted to label our model using $\Sigma=\{p,\nec p\}$, we would have to set $\ell(-1)=\{p,\nec p\}$ and $\ell(0)=\{p\}$. Nevertheless, by the semantics of $\nec p$, we know that $f^n(0)\in \val p$ for all $n<\omega$, yet this information is not recorded in $\ell(0)$.

We could get around this issue by using the classical semantics and the G\"odel-Tarski translation, where $(\nec p)^\blacksquare=\blacksquare\nec \blacksquare p$. If we use $\Sigma={\rm sub}(\blacksquare\nec \blacksquare p)$, we would have $\blacksquare\nec \blacksquare p \in \ell(-1)$, but also $\nec \blacksquare p \in \ell(-1)$ and hence $\nec \blacksquare p \in \ell(0)$. With this, $0$ would `remember' that all of its temporal successors must satisfy $\blacksquare p$. However, note that $\nec \blacksquare p \not \in \ell(x)$ for any $x>0$, and hence $\ell$ is no longer continuous. This is a problem for us, since the continuity of $\ell $ is used in an essential way (for example, in the proof of Lemma \ref{LemmNotTall}) to bound the size of our quasimodels.
\end{exam}

Another useful class of labelled systems is given by quasimodels:

\begin{defn}\label{ndqm}
Given $\Sigma \Subset \landi$, a {\em weak $\Sigma$-quasimodel} is a $\Sigma$-labelled system $\mathfrak Q$ such that $\mathcal T_\mathfrak Q$ is equal to the down-set topology for a partial order which we denote $\peq_\mathfrak Q$. If moreover $\mathfrak Q$ is a well $\Sigma$-labelled system, then we say that $\mathfrak{Q}$ is a {\em $\Sigma$-quasimodel.}
\end{defn}

\begin{exam}\label{ExQuasi}
Recall from Example \ref{ExAlekx} that the formula
\[\varphi=\forall(\ineg p \vee \ps p) \imp \ineg \ps p \vee \ps p \]
is valid on any Alexandroff dynamical model. However, Figure \ref{example} exhibits a quasimodel\footnote{Compare to \cite{me2}, where we use a very similar quasimodel to falsify the formula $\psi=\nec\blacksquare p\to \blacksquare \nec p$, also valid over Alexandroff systems. Note, however, that this $\psi$ uses the classical semantics and thus we use a slightly more complicated example here.} falsifying $\varphi$. Let $\theta = \ineg \ps p \vee \ps p $, and define $\ell$ by
\begin{align*}
\ell(u) & =\{\bm {\forall( \ineg p\vee \ps p)},\bm {( \ineg p\vee \ps p)}, \bm{\ineg p}\}\\
\ell(v) & =\{\bm{ \ps p},\bm{\ineg p}, \bm{ ( \ineg p\vee \ps p ) } , \varphi , \theta ,\forall( \ineg p\vee \ps p)  \}\\
\ell(w) & =\{ \bm{ p},\bm { ( \ineg p\vee \ps p)},\varphi,\theta ,\forall( \ineg p\vee \ps p), \ps p\}.\\
\end{align*}
We include the full labels for completeness, but the most relevant formulas are displayed in boldface. In particular, observe that $p$ belongs only to $\ell(w)$, and $\varphi\not\in \ell(u)$; the latter means that $u$ falsifies $\varphi$, and thus $\varphi$ is not valid over the class of ${\rm sub}(\varphi)$-quasimodels. As we will see, this implies that $\varphi$ is not valid over the class of dynamical systems.
\begin{figure}
\begin{center}

\begin{tikzpicture}[scale=0.7]

\draw[thick] (0,0) circle (.35);

\draw (0,0) node {$u$};

\draw[thick,->] (0,+.5) arc (0:+270:.5) ;
 
\draw (-1.1,+1.1) node {$S$};

\draw[thick] (0,-3) circle (.35);

\draw (0,-3) node {$v$};

\draw[thick,->] (0,-3.5) arc (0:-270:.5) ;
 
\draw (-1.1,-4.1) node {$S$};

\draw[thick,-> ] (.5,-3) -- (2.5,-3);

\draw (1.5,-3.5) node {$S$};

\draw[thick,<-] (0,-2.5) -- (0,-.5);

\draw (-.5,-1.5) node {\large$\rotatebox[origin=c]{90}{$\peq$}$};

\draw[thick] (3,-3) circle (.35);

\draw (3,-3) node {$w$};

\draw[thick,->] (3,-3.5) arc (-180:+90:.5) ;
 
\draw (4.1,-4.1) node {$S$};

\end{tikzpicture}

\end{center}
\caption{A non-deterministic quasimodel.}\label{example}
\end{figure}
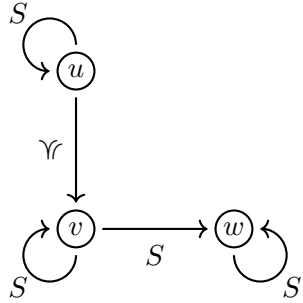
\end{exam}

\section{Producing Dynamic Topological Models from Quasimodels}\label{SecGen}

Note that if $\mathfrak Q$ is a quasimodel, then $S_\mathfrak Q$ is not necessarily a function, and thus we may not view $\mathfrak Q$ directly as a dynamical topological model. However, we can extract a model $\lm$ from it via an unwinding construction, and we call the resulting structure the {\em limit model} of $\mathfrak Q$. The general idea is to consider infinite `deterministic' paths on $\mathfrak Q$ as points in the limit model; however, we will only keep those paths $\vec w$ with the property that, if $\ps\varphi$ occurs in $\vec w$, then $\varphi$ must also occur at a later time. These are the realizing paths of $\mathfrak Q$.

\subsection{Realizing paths}
A {\em path} in a quasimodel $\mathfrak Q$ is any finite sequence $\<w_n\>_{n < \alpha}$ with $\alpha \leq \omega$ such that $w_n \mathrel S w_{n+1}$ whenever $n+1 < \alpha$. An infinite path
${\vec w}=\< w_n\>_{n<\omega}$
is {\em realizing} if for all $n<\omega$ and $\ps\psi\in \ell(w_n)$, there exists $k\geq n$ such that $\psi\in \ell  (w_k)$.

Denote the set of realizing paths by $ \domlm $. This will be the domain of the limit model of $\mathfrak Q$. The transition function on $ \domlm $ will be the `shift' operator, defined by
$\Slm (\<w_n\>_{n<\omega} )=\<w_{n+1}\>_{n<\omega}.$
This simply removes the first element in the sequence.

For our construction to work we must guarantee that there are `enough' realizing paths. The following definition makes this precise:

\begin{defn}
Let $\Sigma\Subset\landi$, $\mathfrak Q$ be a $\Sigma$-quasimodel, and $E\subseteq \domlm$.
Then, $E$ is {\em extensive} if
\begin{enumerate}
\item $E$ is closed under $\Slm$, and
\item\label{ItExtTwo} any finite path
$\<w_0,w_1,..,w_m\>$
in $\mathfrak Q$ can be extended to an (infinite) realizing path
$
{\vec w}=\< w_i\>_{i<\omega}\in E.
$
\end{enumerate}
\end{defn}
\begin{lem}\label{extension}
Let $\Sigma$ be a finite set of formulas closed under subformulas. If
$\mathfrak Q$ is a $\Sigma$-quasimodel, then $\domlm$ is extensive.
\end{lem}

\proof
It is obvious that $\domlm$ is closed under $\Slm$, so we focus on \eqref{ItExtTwo}.
Let
$\<w_i\>_{i\leq m}$
be a finite path and
$(\psi_i)_{i < n}$
be all formulas such that $\ps\psi_i\in \ell(w_m)$ (so that $n=0$ if there are no such $\psi_i$'s). Let $k_0=0$ and suppose inductively that we have constructed a path
\[\<w_0,\ldots,w_{m},\ldots, w_{m+k_i}\>\]
such that, for all $j\in [1,i)$, $\psi_j\in \ell(w_t)$ for some $t\in [m,m+k_i]$. Because $S=S_\mathfrak Q$ is $\omega$-sensible, if $\ps\psi_i\in \ell(w_{m+k_i})$, we know that there exist $k'_{i}\geq 0$ and $v_i\in S^{k'_{i }}(w_{m+k_i})$ such that $\psi_i\in \ell(v_i)$; otherwise, set $k'_{i }=0 $ and $v_{i} = w_{m+k_i}$. We then define $k_{i+1}=k_i+k'_{i}$ and choose
\[w_{m+k_i}\mathrel S w_{m+k_i+1}\mathrel S \hdots \mathrel S  w_{m+k_{i+1}}=v_i.\]
Repeating this process, we obtain a path
\[\<w_0,...,w_m,w_{m+1}\hdots w_{m+k_{n}}\>\]
realizing all $\ps\psi \in \ell(w_m)$. Note that we may have $k_{n}=0$ (for example, if $n=0$), in which case we choose arbitrary $w_{m+1}$ such that $w_m\mathrel S w_{m+1}$, using the seriality of $S$. We then repeat the process beginning with $\<w_0,...,w_{m+k_{n}}\>$, and continue countably many times until an infinite realizing path is formed.
\endproof

If $\mathfrak A$ is an Alexandroff model and $v\peq w$, then for all $i<\omega$ we also have that $f^i_\mathfrak A(v) \peq  f^i_\mathfrak A(w)$, and moreover $ \big ( f^i_\mathfrak A(v) \big)_{i<\omega}$ is a realizing path. However, this does not always hold when we replace $f_\mathfrak A$ by a non-deterministic relation: if $\mathfrak Q$ is a quasimodel, $(w_i)_{i<\omega}$ is a realizing path, and $v_0\peq w_0$, it may be that we cannot `complete' $v_0$ to a realizing path $\vec v = (v_i)_{i<\omega}$ so that $v_i \peq w_i$ for {\em all} $i$. Fortunately, we do not need such a path in our limit model; it is sufficient to `approximate' it by producing paths $\vec v$ so that $v_i \peq w_i$ for all $i < n$, provided $n$ can be taken arbitrarily large. The next lemma will help us do this. We omit the proof, which follows by an easy induction using the continuity of $S_\mathfrak Q$.

\begin{lem}\label{pathcont}
Let
$\mathfrak{Q}$
be a $\Sigma$-quasimodel,
$\<w_i\>_{i \leq n}$
a finite path, and $v_0$ be such that $v_0\peq w_0.$ Then, there exists a path
$\<v_i\>_{i \leq n}$
such that, for $i \leq n$, $v_i\peq w_i$.
\end{lem}

In fact, the path $\vec v$ of Lemma \ref{pathcont} can be used to approximate an infinite path $\vec w$ in a rather precise way; to formalize this, we need to endow $\domlm$ with a topology.

\subsection{Limit models}

If $\Sigma\Subset\landi$ and
$\mathfrak Q$
is a $\Sigma$-quasimodel, the relation ${\peq}={\peq_\mathfrak Q}$ induces a topology on $|\mathfrak Q|$, as we have seen before, by letting open sets be those which are downwards closed under $\peq$. Likewise, $\peq$ induces a very different topology on $\domlm$, in a rather natural way:

\begin{lem}\label{basis}
For each ${\vec w}\in \domlm$ and $n < \omega$,
define
\[\basic n {\vec w}=\Big \{ ( v_i )_{i<\omega}\in \domlm:\forall i\leq n, v_i\acc w_i \Big \} .\]
Then, the set
$\mathcal B_\peq=\big \{ \basic n {\vec w}:{\vec w}\in \domlm, n < \omega \big \}$
forms a topological basis on $\domlm$.
\end{lem}
\proof 
To check that it satisfies \ref{DefBasis}.\ref{ItBasisOne}, note that it is obvious that, given any path ${\vec w}\in \domlm,$ there is a basic set containing it (say, $\basic 0{\vec w}$), and hence
$\domlm=\bigcup_{{\vec w}\in\domlm}\basic 0{\vec w}.$ As for \ref{DefBasis}.\ref{ItBasisTwo}, assuming that ${\vec w}\in \basic{n}{\vec u}\cap \basic{m}{\vec v},$ we observe that
$  { \basic{\max(n,m)}{\vec w} } \subseteq { \basic{n}{\vec u}\cap \basic{m}{\vec v}}$, as needed.\endproof

\begin{defn}
The topology $\Tlm$ on $\domlm$ is the topology generated by the basis $\mathcal B_\peq$.
\end{defn}

Now that we have equipped the set of realizing paths with a topology, we need a continuous transition function on it to have a dynamical system. Our `shift' operator $\Slm$ will do the trick.

\begin{lem}\label{lemSisCont}
The function
$\Slm: \domlm \to \domlm$
is continuous under the topology $\Tlm$.
\end{lem}
\proof Let
${\vec w}=( w_i)_{i<\omega}$
be a realizing path and $\basic n{\Slm ({\vec w})}$ be a neighborhood of $\Slm({\vec w})$. Then, if
${\vec v}\in \basic{n+1}{\vec w},$
$w_i\acc v_i$ for all $i\leq n+1$, so $w_{i+1}\acc v_{i+1}$ for all $i\leq n$ and
$\Slm(\vec v)\in \basic n{\Slm(\vec w)}.$
Hence
$\Slm (\basic{n+1}{\vec w})\subseteq \basic n{\Slm(\vec w)},$
and $\Slm$ is continuous.\endproof

Finally, we will use $\ell$ to define a valuation: if $p$ is a propositional variable, set
\[\val p^\ell= \big \{ {\vec w}\in \domlm : p \in \ell(w_0) \big \}.\]
With this, we are now ready to assign a dynamic topological model to each quasimodel:

\begin{defn}
Given a $\Sigma$-quasimodel $\mathfrak{Q},$ define
\[\vec {\mathfrak{Q}}=\Big ( \domlm,\Tlm,\Slm,\lb\cdot\rb^\ell \Big ) \]
to be the {\em limit model} of $\mathfrak{Q}$.
\end{defn}

Of course this model is only useful if $\lb\cdot\rb^\ell$ matches with $\ell$ on all formulas of $\Sigma$, not just propositional variables. Fortunately, this turns out to be the case.

\begin{lem}\label{sound}
Let $\Sigma\Subset\landi$, $\mathfrak Q$ be a $\Sigma$-quasimodel,
${\vec w}=( w_n)_{n < \omega }\in \domlm$
and $\varphi\in \Sigma$.
Then,
\[\lb\varphi\rb^\ell=\{\vec w : \varphi\in \ell(w_0)\}.\]
\end{lem}
\proof The proof goes by standard induction of formulas. The induction steps for $\wedge,\vee $ are immediate; here we will only treat the cases for $\imp,\nx,\ps,\nec$.

\paragraph{Case 1: $\varphi=\psi\imp\theta$}
If $\varphi\in \ell(w_0)$, take the neighborhood $\basic 0 {w_0}$ of ${\vec w}$ and consider $\vec v\in \basic 0{w_0}$, so that $v_0\acc w_0$. Since $\mathfrak Q$ is a $\Sigma$-labelled space, it follows that, if $\psi\in \ell(v_0)$, then $\theta\in \ell(v_0)$; by the induction hypothesis, this means that if $\vec v\in \lb\psi\rb^\ell$, then $\vec v\in \lb\theta\rb^\ell$. Since $\vec v\in \basic 0{w_0}$ was arbitrary, we conclude that $\vec w\in \lb\psi\imp\theta\rb^\ell$.

On the other hand, if $\psi\imp \theta\not\in \ell(w_0)$, there is $v_0\peq w_0$ such that $\psi\in \ell(v_0)$ but $\theta\not\in \ell (v_0)$. Consider any basic neighborhood $\basic{n}{\vec w}$ of $\vec w$. Then, by Lemma \ref{pathcont}, there exists a path
$\<v_0,...,v_n\>\subseteq |\mathfrak Q|$
such that, for all $i\leq n$, $v_i\acc w_i$ and for $i<n$, $v_i\mathrel S v_{i+1}$. Because $\domlm$ is extensive,
$ ( v_i ) _{i\leq n}$
can be extended to a realizing path ${\vec v}\in \domlm$. Then, ${\vec v}\in \basic n{{\vec w}}$, and by induction hypothesis 
$\vec v\in \lb \varphi\rb^\ell\setminus \lb \theta\rb^\ell.$ Since $n$ was arbitrary, we conclude that $\vec w\not\in \lb\psi\imp\theta\rb^\ell.$

\paragraph{Case 2: $\varphi=\nx\psi$} This case follows from the fact that $(w_0,w_1)$ is sensible and the induction hypothesis.

\paragraph{Case 3: $\varphi=\ps\psi$} 
Because ${{\vec w}}$ is a realizing path, if $\ps\psi\in \ell(w_0)$, then $\psi\in \ell (w_n)$ for some $n\geq 0$. We can use the induction hypothesis to conclude that
${\Slm}^n ( \vec w )\in\val\psi^\ell$,
and so ${{\vec w}}\in\val{\ps\psi}^\ell$.

For the other direction, assume that $\ps\psi\not\in \ell(w_0)$. For all $n$, $(w_n,w_{n+1})$ is sensible so an easy induction shows that $\psi\not\in \ell(w_n)$, so that by the induction hypothesis ${\Slm}^n({\vec w})\not\in \val\psi^\ell;$
since $n$ was arbitrary, ${\vec w}\not\in \val{\ps\psi}^\ell.$

\paragraph{Case 4: $\varphi=\forall\psi$} If $\forall \psi\in \ell(w_0)$, then since $\ell$ is honest, we can use the induction hypothesis to conclude that $\vec v\in\val\psi^\ell$ for all $\vec v\in \domlm$, and thus $\vec w\in\val{\forall\psi}^\ell$. If $\forall\psi\not\in \ell(w_0)$, there is $v_0\in|{\mathfrak Q} |$ such that $\psi\not\in \ell(v_0)$. We can extend $v_0$ to a realizing path $\vec v$ and, by the induction hypothesis, $\vec v\not\in \val\psi^\ell$, hence $\vec w\not\in \val{\forall\psi}^\ell$.
\endproof

We are now ready to prove the main result of this section, which in particular implies that $\icltl$ is sound for the class of quasimodels.

\begin{prop}\label{second}
Let $\Sigma\Subset\landi$ and $\mathfrak Q$ be any quasimodel. Then, $\lm$ is a model, and for all $\varphi \in \Sigma$ we have that $\lm$ falsifies $\varphi$ if and only if $\mathfrak Q$ falsifies $\varphi$.
\end{prop}

\proof
It is clear from its definition that $\Slm$ is a function and by Lemma \ref{lemSisCont} it is continuous, so that $(\domlm,\Tlm,\Slm)$ is a dynamical system.
Now, suppose that $w_0\in |\mathfrak Q|$ is such that $\varphi \in \Sigma \setminus \ell(w_0)$. By Lemma \ref{extension}, $w_0$ can be extended to a realizing path ${\vec w}$. It follows from Lemma \ref{sound} that ${\vec w}\not\in \val\varphi^\ell,$ and therefore $\varphi$ is falsified on $\lm$. Conversely, if $\vec w = (w_i)_{i<\omega} \in \domlm$ falsifies $\varphi$ then $\varphi \not \in \ell(w_0)$, so that $\fr Q$ falsifies $\varphi$ if and only if $\lm$ does.
\endproof

\begin{exam}\label{ExQ}
Recall from Examples \ref{ExAlekx} and \ref{ExQuasi} that the formula
\begin{equation}\label{EqAleksForm}
\varphi=\forall(\ineg p \vee \ps p) \imp \ineg \ps p \vee \ps p
\end{equation}
is valid over the class of Alexandroff systems, but is falsifiable in a finite quasimodel. Thus it follows from Proposition \ref{second} that $\varphi$ is also falsifiable in a dynamical topological model.

It is instructive to construct such a model directly, rather than appealing to the theorem. To this end, let us refute $\varphi$ on $\mathbb R$. Define $f\colon \mathbb R \to \mathbb R$ by $f(x) = 2x$, and $\val p = (1,\infty)$. Then it is not hard to check that $\val{\ineg p} = (-\infty, 1)$ and $\val {\ps p} = (0,\infty)$, so that $\val{ \forall (\ineg p\vee \ps p)} = \mathbb R$. On the other hand $\val{\ineg \ps p} = (-\infty, 0)$, so that $0 \not \in \val{\ineg \ps p \vee \ps p}$ and $\varphi$ is not valid on $\mathbb R$.
\end{exam}

In view of Example \ref{ExAlekx}, we conclude that $\icltl\subsetneq{\itle}$. Note however that not every falsifiable formula is falsifiable on $\mathbb R$, given that the latter is connected, in the following sense:

\begin{defn}
A topological space $\mathfrak X$ is {\em disconnected} if there are non-empty, disjoint open sets $A,B\subseteq |\mathfrak X|$ such that $A\cup B = |\mathfrak X|$; otherwise, $\mathfrak X$ is {\em connected.}
\end{defn}

It is well-known that $\mathbb R^n$ is connected for any $n<\omega$. As noted by Shehtman, the formula
\begin{equation}\label{EqConn}
\forall(p\vee\ineg p)\imp (\forall p\vee\forall \ineg p)
\end{equation}
is valid over the class of connected spaces \cite{ShehtmanEverywhereHere}, from which the validity of $\varphi$ for systems based on a connected space is a straightforward consequence.

\section{Moments}\label{SecFour}

We have seen that any valid $\varphi\in\Sigma$ is valid over the class of $\Sigma$-quasimodels. The converse is true, and for this we need to construct, given a model $\mathfrak X$ falsifying $\varphi$, a quasimodel also falsifying $\varphi$. This quasimodel will be denoted $\mathfrak X/\Sigma$. The worlds of $\mathfrak X/\Sigma$ will be called $\Sigma$-moments; the intuition is that a $\Sigma$-moment represents all the information that holds at the same `moment of time'.

\begin{defn}
A {\em $\Sigma$-moment} is a $\Sigma$-labelled space $\mathfrak w$ such that $\mathcal T_\mathfrak w$ is the down-set topology of a partial order $\peq_\mathfrak w$, and $( |\mathfrak w|,\peq_\mathfrak w )$ is a finite tree with (unique) root $\Root{\mathfrak w}$. We will write $\ell_\Sigma(\mathfrak w)$ instead of $\ell_{\mathfrak w}(\Root{\mathfrak w})$. The set of $\Sigma$-moments is denoted $M_\Sigma$.
\end{defn}

\begin{defn}
Let $\mathfrak w$ be a $\Sigma$-moment. For $w\in |\mathfrak w|$, let $\mathfrak w [ w ] = {\fw\upharpoonright \mathord\downarrow w}$, i.e.,
\[\mathfrak w [ w ] = \big ( \ \mathord\downarrow w \ , \ {\peq_\mathfrak w} \upharpoonright \mathord\downarrow w  \ , \ \ell_\mathfrak w   \upharpoonright \mathord\downarrow w \ \big ) .\]
We write $\mathfrak w\peq_\Sigma \mathfrak v$ if $\mathfrak v=\mathfrak w [ w ] $ for some $w\in |\mathfrak w|$.
\end{defn} 

We now wish to define a weak $\Sigma$-quasimodel $\moments\Sigma$ over $M_\Sigma$. For this, it remains to define a sensible relation on $M_\Sigma$.

\begin{defn}\label{ts}
Say $\mathfrak w$ is a {\em temporal successor} of $\mathfrak v$, denoted $\mathfrak v \mathrel S_\Sigma \mathfrak w$, if there exists a sensible relation
$R\subseteq |\mathfrak v|\times |\mathfrak w|$
such that $\Root{\mathfrak v}\mathrel R \Root{\mathfrak w}$.
\end{defn}

\begin{lem}\label{LemmIsSensible}
Let $\Sigma\Subset\landi$ and $\fw,\fw',\fv\in M_\Sigma$.
\begin{enumerate}

\item If $\fw\circrel\fv$ then $(\ell_\Sigma(\fw),\ell_\Sigma(\fv))$ is sensible.

\item \label{ItIsSensibleTwo}

If $\fw'\peq_\Sigma\fw\circrel \fv$, there is $\fv'\peq_\Sigma \fv$ such that $\fw'\circrel\fv'$.

\end{enumerate}

\end{lem}

\proof
Suppose that $\fw\circrel\fv$, and let $R\subseteq |\mathfrak v|\times |\mathfrak w|$ be a sensible relation
such that $\Root{\mathfrak v}\mathrel R \Root{\mathfrak w}$. Then, the pair $(\ell_\fw(\Root \fw),\ell_\fv(\Root \fv))$ is sensible because $R$ is sensible, but it is equal to $(\ell _\Sigma(  \fw),\ell_\Sigma (  \fv))$, as needed for the first claim.

For the second, suppose further that $\fw'\peq_\Sigma \fw$. This means that $\fw'=\fw[w]$ for some $w\in |\fw|$. Since $R$ is sensible and $|\fv|$ is open, $R^{-1} |\fv|$ is open as well, meaning in particular that $w\mathrel R v$ for some $v\in |\fv|$. Now consider $R'=R\cap (\mathop\downarrow w\times \mathop\downarrow v)$. Since $\mathop\downarrow v$ is open it follows that $R'$ is continuous as well, and every pair in $R'$ is sensible since every pair in $R$ was. Moreover, $\Root{\fw[w]}=w\mathrel R' v= \Root{\fv[v]}$. Thus $R'$ witnesses that $\fw'=\fw[w]\circrel \fv[v]\peq_\Sigma\fv$, as needed.
\endproof

Lemma \ref{LemmIsSensible}.\ref{ItIsSensibleTwo} essentially says that $S_\Sigma$ is a continuous relation with respect to the down-set topology induced by $\peq_\Sigma$ (see Figure \ref{FigContRel}). 
With this, we are ready to define our `canonical' weak quasimodels.

\begin{defn}
Let $\Sigma\Subset\landi$. Define $M_\Sigma$ to be the set of all finite $\Sigma$-moments, and set
$\moments\Sigma=\<M_\Sigma,{\peq_\Sigma},S_\Sigma,\ell_\Sigma\>$.
\end{defn}

\begin{lem}
If $\Sigma\Subset\landi$, then $\moments\Sigma$ is a weak quasimodel.
\end{lem}

\proof That $\circrel$ is sensible is Lemma \ref{LemmIsSensible}, so it remains to show that $(M_\Sigma,{\peq_\Sigma},\ell_\Sigma)$ is a $\Sigma$-labelled space. First we check that $\ell_\Sigma$ is continuous. This amounts to showing that, if $\fw\cca_\Sigma\fv\in M_\Sigma$, then $\ell_\Sigma(\fw)\subseteq \ell_\Sigma(\fv)$. But, $\fw\cca_\Sigma\fv$ means that $\fv=\fw[v]$ for some $v\peq_\fw \Root\fw$, hence $\ell_\Sigma(\fv)=\ell_\fw (v)\supseteq \ell_\fw(\Root \fw)=\ell_\Sigma(\fw)$ by the continuity of $\ell_\fw$. Similarly, if $\delta = ( \varphi \imp \psi ) \in \defect{\ell_\Sigma(\fw)}$, then there is $v\peq_\fw \Root\fw$ such that $v$ revokes $\delta$; but then, $\fw[v]\peq_\Sigma \fw$ and $\fw[v]$ revokes $\delta$.
\endproof

Henceforth, we may write $\peq,S,\ell$ instead of $\peq_\Sigma,S_\Sigma,\ell_\Sigma$ when this does not lead to confusion. Observe that $S_\Sigma$ is not necessarily $\omega$-sensible, and $\ell_\Sigma$ is not necessarily honest, so $\moments\Sigma$ is not a quasimodel as it stands. It does, however, contain substructures which are proper quasimodels, as we will see later.

\subsection{Building moments from smaller moments}

Often we will want to construct a $\Sigma$-moment from smaller moments. Here we will define the basic operation we will use to do this, and establish the conditions that the pieces must satisfy. Below, $\coprod$ denotes a disjoint union.

\begin{defn}
Let $\varphi$ be a formula of $\icltl$, $\Sigma\Subset\landi$, $\Phi\in\type \Sigma$ and $U\subseteq \moments\Sigma$.

Define
$\mathfrak w=\add\Phi U$
by setting
\begin{align*}
 |\mathfrak w|&=\{\Phi\}\cup \coprod_{\mathfrak u\in U}|{\mathfrak u}|,\\
\acc_\mathfrak w&=\big (\{\Phi\}\times |\mathfrak w|\big )\cup \coprod_{\mathfrak u\in U}\acc_{\mathfrak u},\\
\ell_\mathfrak w(v)&=
\{(\Phi,\Phi)\}\cup \coprod _{\mathfrak u\in U}\ell_{\mathfrak u}(v).
\end{align*}

\end{defn}

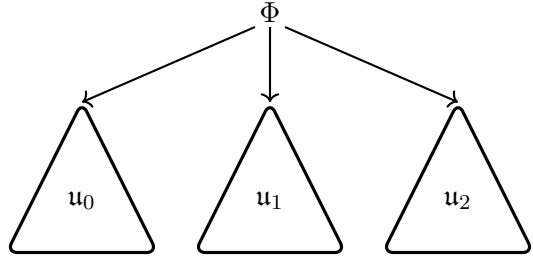
\begin{figure}
\begin{center}

\begin{tikzpicture}

\draw (2.5,1.2) node {$\Phi$};

\draw[thick,->] (2.3,1) --  (0,0);

\draw[rounded corners, very thick] (0,0) -- (-1,-2) -- (1,-2) -- cycle;

\draw (0,-1.3) node {$\fu_0$};

\draw[thick,->] (2.5,+1) --  (2.5,0);

\draw[rounded corners, very thick] (2.5,0) -- (1.5,-2) -- (3.5,-2) -- cycle;

\draw (2.5,-1.3) node {$\fu_1$};

\draw[thick,->] (2.7,+1) --  (5,0);

\draw[rounded corners, very thick] (5,0) -- (4,-2) -- (6,-2) -- cycle;

\draw (5,-1.3) node {$\fu_2$};

\end{tikzpicture}

\end{center}
\caption
{
The moment $\add\Phi U$, where $U=\{\fu_0,\fu_1,\fu_2\}$ and $\Phi\subseteq \Sigma\Subset\landi$.
}
\end{figure}

Note that $\add\Phi U$ will not always be a moment, since $\ell_{\add\Phi U}$ thus defined might not be continuous, or it might not revoke some defect. To ensure that we do obtain a new moment, we need for $(\Phi,U)$ to be a kit.

\begin{defn}\label{cohe}
Fix $\Sigma\Subset\landi$.
\begin{enumerate}

\item

A {\em $\Sigma$-kit} is a pair $(\Phi, U)$, with $\Phi\in \type\Sigma$ and $U\subseteq M_\Sigma$ finite, such that
\begin{enumerate}

\item $\Phi\subseteq \ell (\mathfrak u)$ for all $\mathfrak u\in U$, and 

\item whenever $\varphi\imp\psi\in \defect \Phi$, there is $\mathfrak u\in U$ such that $\varphi\imp\psi\not\in \ell (\mathfrak u)$.

\end{enumerate}

\item

If $\mathfrak w\in\moments\Sigma$, we say that the $\Sigma$-kit $(\Phi,U)$ is a {\em successor kit for $\mathfrak w$} if
\begin{enumerate}

\item the pair $(\ell(\mathfrak w),\Phi)$ is sensible, and

\item if $\fv\prec\fw$, there is $\mathfrak u\in U$ such that $\mathfrak v\circrel\mathfrak u$.

\end{enumerate}
\end{enumerate}
\end{defn}

\begin{lem}\label{add}
Let $\Sigma\Subset\landi$, $\Phi\subseteq \Sigma$, and $V\subseteq \moments \Sigma$.

Then,
\begin{enumerate}
\item\label{ItAddOne} $\add \Phi V$ is a $\Sigma$-moment if and only if $(\Phi,V)$ is a $\Sigma$-kit, and\vskip 5pt
\item\label{ItAddTwo} $\add \Phi V$ is a $\Sigma$-moment with $\mathfrak w\circrel\add \Phi V$ whenever $(\Phi, V)$ is a successor kit for $\mathfrak w$.
\end{enumerate}
\end{lem}

\proof
To prove \eqref{ItAddOne}, one can check that the conditions for being a $\Sigma$-kit correspond exactly to the conditions in Definition \ref{frame}. For \eqref{ItAddTwo}, let $w_1,\hdots,w_n$ enumerate the set $|\fw|\setminus\{\Root{\mathfrak w}\}$, and for $i\leq n$, let $\mathfrak w_i=\mathfrak w[ w_i]$. For $i\leq n$, let $\mathfrak v_i\in V$ be such that $\mathfrak w_i\circrel\mathfrak v_i$, and $R_i\subseteq |\mathfrak w_i|\times |\mathfrak v_i|$ be sensible. Then, define $R\subseteq |\mathfrak v|\times \left |\add \Phi V\right |$
by $w\mathrel R v$ if and only if
\begin{enumerate*}[label=(\alph*)]

\item $w=\Root{\mathfrak w}$ and $v=\Root{\mathfrak v}$, or

\item $w\in |\mathfrak w_i|$, $v\in |\mathfrak v_i|$ and $w\mathrel R_i v$.

\end{enumerate*}
It is readily seen that the relation $R$ thus defined is sensible.
\endproof

\section{Simulations}\label{SecSim}

In order to prove that every falsifiable formula is falsifiable in a quasimodel, we need to represent dynamical models using quasimodels. Simulations are relations between quasimodels and models (or labelled spaces in general), and provide the basic tools we need to achieve such a representation.

\begin{defn}
Let $\Sigma\subseteq \Delta \subseteq \landi$ both be closed under subformulas, $\mathfrak X$ be a $\Sigma$-labelled space and $\mathfrak Y$ be a $\Delta$-labelled space. A continuous relation
$\chi\subseteq |\mathfrak X|\times |\mathfrak Y|$
is a {\em simulation} if, for all $(x,y)\in \chi$, $\ell_\mathfrak X (x)= \ell_\mathfrak Y(y) \cap \Sigma.$
\end{defn}

We will call the latter property {\em label-preservation}. Note that it may be that $\ell_\mathfrak X$ is honest while $\ell_\mathfrak Y$ is not, or vice-versa. However, this cannot happen with a total, surjective simulation:

\begin{lem}\label{LemTotSim}
Suppose that $\Sigma\subseteq \Delta \subseteq \landi$ are both closed under subformulas, $\mathfrak X$ is a $\Sigma$-labelled space, $\mathfrak Y$ is a $\Delta$-labelled space, and $\chi\subseteq |\mathfrak X|\times |\mathfrak Y|$
is a total, surjective simulation. Then, if $\ell_\mathfrak Y$ is honest, it follows that $\ell_\mathfrak X$ is honest as well.
\end{lem}

We omit the proof, which is straightforward; note that if $\Sigma = \Delta$ then the last line becomes an equivalence. Simulations are useful for comparing the purely topological behavior of labelled structures. However, to compare their dynamic behavior, we need something a bit stronger.

\begin{defn}\label{simulation}
Let $\mathfrak X$ and $\mathfrak Y$ be labelled systems. A {\em dynamic simulation} between $\mathfrak X$ and $\mathfrak Y$ is a simulation $\chi\subseteq |\mathfrak X|\times |\mathfrak Y|$ such that $
S_\mathfrak Y\chi\subseteq\chi S_\mathfrak X.
$
\end{defn}

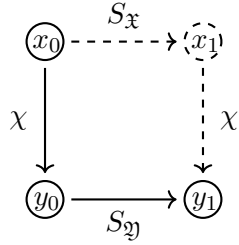
\begin{figure}

\begin{center}

\begin{tikzpicture}[scale = 0.7]

\draw[thick]  (0,0) circle (.35);

\draw (.04,-.05) node {$y_0$};

\draw[thick, ->] (.5,0) -- (2.5,0);

\draw (1.5,-.5) node {$S_\mathfrak Y$};

\draw[thick] (3,0) circle (.35);

\draw (3+.04,-.05) node {$y_1$};

\draw[thick] (0,3) circle (.35);

\draw (.04,3-.04) node {$x_0$};

\draw[thick, ->,dashed] (.5,3) -- (2.5,3);

\draw (1.5,3.5) node {$S_\mathfrak X$};

\draw[thick, ->] (0,2.5) -- (0,.5);

\draw (-.5,1.5) node {$\chi$};

\draw[thick, dashed] (3,3) circle (.35);

\draw (3+.05,3-.04) node {$x_1$};

\draw[thick,->,dashed] (3,2.5) -- (3,.5);

\draw (3.5,1.5) node {$\chi$};

\end{tikzpicture}

\end{center}

\caption{The above diagram can always be completed if $\chi \subseteq |\mathfrak X| \times |\mathfrak Y|$ is a dynamic simulation.}\label{figSim}

\end{figure}

The following properties are readily verified:

\begin{lem}\label{LemmPropSim}

Let $\mathfrak X,\mathfrak Y,\mathfrak Z$ be labelled systems and $\chi\subseteq |\mathfrak X|\times  |\mathfrak Y|$, $\xi\subseteq |\mathfrak Y|\times  |\mathfrak Z|$ be simulations. Then:

\begin{enumerate}

\item\label{ItPropSimComp} $\xi\chi\subseteq |\mathfrak X|\times  |\mathfrak Z|$ is a simulation. Moreover, if both $\chi$ and $\xi$ are dynamic, then so is $\xi\chi$.

\item \label{ItPropSinSub} If $U\subseteq |\mathfrak X|$ and $V\subseteq |\mathfrak Y|$ are open, then $\chi\upharpoonright U\times V$ is a simulation.

\item\label{ItPropSimUn} If $\Xi\subseteq\mathcal P(|\mathfrak X|\times  |\mathfrak Y|)$ is a set of simulations, then $\bigcup \Xi$ is also a simulation.

\end{enumerate}

\end{lem}

If $\chi$ is a dynamic simulation between a weak quasimodel $\mathfrak W$ and a dynamical model $\mathfrak X$, it does not immediately follow that $S_\mathfrak W$ is $\omega$-sensible. However, we can use the simulation $\chi$ to extract an $\omega$-sensible substructure from $\mathfrak W$.

\begin{lem}\label{isquasi}
Let $\mathfrak W$, $\mathfrak V$ be labelled systems and $\chi\subseteq |\mathfrak W|\times |\mathfrak V|$ be a dynamical simulation. Then, if $S_\mathfrak V$ is $\omega$-sensible, it follows that $S_\mathfrak W\upharpoonright \chi^{-1}(|\mathfrak V|)$ is also $\omega$-sensible.
\end{lem}
\proof Let $w\in \chi^{-1}(|\mathfrak V|)$ and $\ps\varphi\in \ell_\mathfrak W(w)$. Since $w\in \chi^{-1}(|\mathfrak V|)$, we can choose $v_0 \in |\mathfrak V|$ such that $w\mathrel \chi v_0$. Since $\ps\varphi\in \ell_\mathfrak W(w)$, it follows that $\ps\varphi\in \ell_\mathfrak V(v)$. Using the assumption that $S_\mathfrak V$ is $\omega$-sensible, we may choose a sequence
\[v_0\mathrel S_\mathfrak V v_1\mathrel S_\mathfrak V v_2\mathrel S_\mathfrak V \hdots \mathrel S_\mathfrak V v_n\]
such that $\varphi\in \ell_\mathfrak V(v_n)$. Since $\chi$ is a dynamic simulation, we may recursively find
\[w=w_0\mathrel S_\mathfrak W w_1\mathrel S_\mathfrak W w_2\mathrel S_\mathfrak W \hdots \mathrel S_\mathfrak W w_n\]
such that $w_i\mathrel \chi v_i$ for each $i\leq n$. Then $\varphi\in \ell_\mathfrak W(w_n)$, as needed.

The seriality of $S_\mathfrak W\upharpoonright \chi^{-1}(|\mathfrak V|)$ is verified similarly, using the seriality of $S_\mathfrak V$.
\endproof

Suppose in particular that $\mathfrak X$ is a model and $x_\ast\in |\mathfrak X|$ falsifies $\varphi\in \Sigma$. Then, if $\mathfrak Q$ is a weak $\Sigma$-quasimodel and $\chi\subseteq |\mathfrak Q|\times |\mathfrak X|$ is a surjective dynamical simulation such that $w_\ast\mathrel \chi x_\ast$ for some $w_\ast\in|\mathfrak Q|$, it follows $\mathfrak Q\upharpoonright \chi^{-1}(|\mathfrak X|)$ is a quasimodel falsifying $\varphi$. Thus our strategy will be to show that there is a weak $\Sigma$-quasimodel $\mathfrak Q$ such that, given {\em any} dynamical model $\mathfrak X$, there is a surjective dynamical simulation $\chi\subseteq |\mathfrak Q|\times |\mathfrak X|$. In principle we could take $\mathfrak Q=\moments\Sigma$, but since we wish to obtain finite quasimodels, it will be convenient to consider a finite substructure of $\moments\Sigma$. The elements of this structure will be the irreducible $\Sigma$-moments, as defined in the next section.

\section{Irreducible Moments}\label{SecIrred}

In order to obtain finite quasimodels, we will restrict $\moments\Sigma$ to moments that are, in a sense, no bigger than they need to be. To be precise, we want them to be minimal with respect to $\redu$, which we define below, along with some other useful relations between moments.

\begin{defn}
Let $\Sigma\Subset\landi$ and $\mathfrak w,\mathfrak v$ be $\Sigma$-moments. We write 

\begin{enumerate}

\item $\mathfrak w\cong\mathfrak v$ if $\mathfrak w,\mathfrak v$ are isomorphic;

\item $\mathfrak w\sqsubseteq \mathfrak v$ if $|\mathfrak w|\subseteq |\mathfrak v|$, ${\peq_\mathfrak w}={\peq_\mathfrak v\upharpoonright |\mathfrak w|}$, and $\ell_\mathfrak w =\ell_\mathfrak v\upharpoonright |\mathfrak w|$;

\item $\mathfrak w\redu\mathfrak v$ if $\mathfrak w\sqsubseteq\mathfrak v$ and there is a continuous, surjective function $\pi\colon |\mathfrak v|\to |\mathfrak w|$ such that $\ell_\mathfrak v(v)=\ell_\mathfrak w(\pi(v))$ for all $v\in |\mathfrak v|$ and $\pi^2=\pi$. We say that $\mathfrak w$ is a {\em reduct} of $\mathfrak v$ and $\pi$ is a {\em reduction}.

\end{enumerate}

\end{defn} 

Note that the condition $\pi^2=\pi$ is equivalent to requiring $\pi(w)=w$ whenever $w\in |\fw|$.

\begin{lem}\label{LemPropRed}

Let $\Sigma\Subset\landi$ and $\fw,\fw',\fv,\fv',\fu\in M_\Sigma$.
\begin{enumerate}

\item\label{ItLemPropRedType} If $\mathfrak w\redu\mathfrak v$, then $\ell(\fw)=\ell(\fv)$.

\item\label{ItLemPropRedSub} If $\fw\redu\fv$ and $w\in |\fw|$, then $\fw[w]\redu \fv[w]$.

\item\label{ItLemPropRedTrans} If $\fw\redu\fv\redu\fu$, then $\fw\redu\fu$.

\item\label{ItLemPropRedNext} If $\fw\circrel\fv$, $\fw'\redu\fw$ and $\fv'\redu\fv$, then $\fw'\circrel\fv'$.

\end{enumerate}

\end{lem}

\proof\leavevmode

\begin{enumerate}

\item[\eqref{ItLemPropRedType}] Let $\pi\colon |\fv|\to |\fw|$ be a reduction. 
Observe that $\Root\fw=\pi(\Root\fv)$, since $\Root\fw \peq \Root \fv$ implies that $\Root\fw=\pi(\Root \fw)\peq \pi(\Root\fv)$, but $\pi(\Root \fv ) \peq \Root\fw$ since it is the root so that the two are equal. It follows that $\ell(\fw)=\ell_\fw(\pi(\Root\fv))=\ell_\fv(\Root\fv)=\ell(\fv)$.

\item[\eqref{ItLemPropRedSub}] Let $\pi\colon |\fv|\to |\fw|$ be a reduction. Then, $ (\pi\upharpoonright \mathord\downarrow w) \colon \mathord \downarrow w \to \mathord\downarrow w$ is label-preserving, continuous and idempotent. It follows that $\fw[w]\redu\fv[w]$.

\item[\eqref{ItLemPropRedTrans}] Let $\pi\colon |\fv|\to |\fw|$ and $\rho\colon |\fu|\to |\fv|$ be reductions. Then, $\pi \rho$ is continuous, surjective, label-preserving, and if $w\in |\fw|$, $\pi \rho (w)=\pi (w)= w$.

\item[\eqref{ItLemPropRedNext}] Let $R\subseteq |\mathfrak w|\times|\mathfrak v|$ be sensible, $\iota\colon |\mathfrak w'|\to|\mathfrak w|$ be the inclusion map, and $\pi\colon|\mathfrak v|\to|\mathfrak v'|$ be a reduction. Then, one readily checks that $\pi R \iota \subseteq |\mathfrak w'|\times |\mathfrak v'|$ is sensible.\qedhere
\end{enumerate}

\endproof

Irreducible moments are the minimal moments under $\redu$. Let us define them formally and establish their basic properties. Below, recall that $\mathfrak v\prec^1\mathfrak w$ means that $\mathfrak v\prec\mathfrak w$, and there is no $\fu$ such that $\fv\prec\fu\prec\fw$.

\begin{defn}
A moment $\fw$ is {\em irreducible} if whenever $\fw \redu\fv$, it follows that $\fw=\fv$.
\end{defn}

\begin{lem}\label{LemExistIrr}
Let $\Sigma\Subset\landi$.
\begin{enumerate}

\item\label{LemExistIrrOne} If $\fw$ is any $\Sigma$-moment, there is $\fv\redu\fw$ such that $\fv$ is irreducible.

\item\label{LemExistIrrTwo} If $\fw$ is irreducible and $\fv\peq\fw$, then $\fv$ is irreducible.

\item\label{LemExistIrrThree} If $\mathfrak w$ is irreducible, $  u,  v\prec^1 r_{\mathfrak w}$ and $\mathfrak w[u] \cong\mathfrak w[v]$, then $  u= v$.

\end{enumerate}
\end{lem}

\proof\leavevmode
\begin{enumerate}

\item[\eqref{LemExistIrrOne}] Choose $\fv\redu\fw$ such that $\# |\fv|$ is minimal. Then, if $\fu\redu\fv$, it follows that $\# |\fu|\leq \# |\fv|$ and, by Lemma \ref{LemPropRed}.\ref{ItLemPropRedTrans}, $\fu\redu\fw$; hence by minimality, $\# |\fu|= \# |\fv|$, and thus $\fu=\fv$. Since $\fu$ was arbitrary, $\fv$ is irreducible.

\item[\eqref{LemExistIrrTwo}] If $\fv\peq\fw$, then $\fv=\fw[w]$ for some $w\in |\fw|$. Let $\pi\colon \mathord\downarrow \to \mathord\downarrow  w$ be a reduction, and define $\rho\colon |\fw|\to |\fw|$ by $\rho(v)=\pi(v)$ if $v\peq w$, $\rho(v)=v$ otherwise. It is readily checked that $\rho$ is a reduction, and since $\fw$ is irreducible, it must be surjective. However, this is only possible if $\pi$ was already surjective.

\item[\eqref{LemExistIrrThree}] Assuming otherwise, let $h\colon |\mathfrak w[u]|\to|\mathfrak w[v]|$ be an isomorphism; then, define $\pi\colon |\mathfrak w|\to|\mathfrak w|$ by $\pi(x)=h(x)$ if $x\in |\mathfrak w[u]|$, $\pi(x)=x$ otherwise. It is readily checked that $\pi$ is a reduction onto a strictly smaller moment.\qedhere

\end{enumerate}
\endproof

The set of irreducible moments forms a weak quasimodel, much like the set of moments $\moments\Sigma$.

\begin{defn}
Let $\Sigma\Subset\landi$. Define $I_\Sigma$ to be the set of all irreducible $\Sigma$-moments, and set $\mathbb I_\Sigma=\moments\Sigma\upharpoonright I_\Sigma$.
\end{defn}

Let us write $\fw\redu_0\fv$ if $\fw\redu\fv$ and $\fv$ is irreducible. Then we have:

\begin{lem}\label{LemRedIsSim}
Let $\Sigma\Subset\landi$. Then,
\begin{enumerate}

\item $\irr\Sigma$ is a weak $\Sigma$-quasimodel, and

\item ${\redu_0}\subseteq |\irr\Sigma|\times |\moments \Sigma|$ is a surjective, dynamic simulation.

\end{enumerate}
\end{lem}

\proof
It is immediate from Lemma \ref{LemExistIrr}.\ref{LemExistIrrTwo} that $I_\Sigma$ is open, and thus $\irr\Sigma$ is a weak quasimodel by Lemma \ref{LemmOpenSubst}. Surjectivity of $\redu_0$ is Lemma \ref{LemExistIrr}.\ref{LemExistIrrOne}, label-preservation is Lemma \ref{LemPropRed}.\ref{ItLemPropRedType}, continuity is Lemma \ref{LemPropRed}.\ref{ItLemPropRedSub}, and dynamicity is Lemma \ref{LemPropRed}.\ref{ItLemPropRedNext}.
\endproof

We wish to show that $\mathbb I_\Sigma$ is always finite. For this, we will first show that irreducible frames cannot be too `tall', which will be a consequence of the next lemma.

\begin{lem}\label{LemmNotTall}
Let $\fw$ be an irreducible $\Sigma$-moment and $v,w\in|\fw|$. Then, if $v\peq w$ and $\ell(v)=\ell(w)$, it follows that $v=w$.
\end{lem}

\proof
 We proceed by induction on ${\rm hgt}(\fw)$: if $w\not=\Root\fw$, then by Lemma \ref{LemExistIrr}.\ref{LemExistIrrTwo}, $\fw[w]$ is irreducible and we can apply the induction hypothesis. So, we may assume otherwise. We may also assume that $v\prec^1 \Root\fw$; for, if $\Root\fw=v$, there is nothing to prove, and if not, choose $v'$ such that $\Root\fw\succ^1 v' \succcurlyeq v$; then, $\ell(\Root\fw)\subseteq \ell(v')\subseteq \ell(v)=\ell(\Root\fw)$, so we also have $\ell(v')=\ell(\Root\fw)$. Then, define $\pi$ by $\pi(v)=\Root\fw$, $\pi(u)=u$ otherwise. It is readily seen that $\pi$ is a reduction and $v\not\in \pi|\fw|$, contradicting the irreducibility of $\fw$.
\endproof

Let us define the {\em height} of $\fw$ as the largest $n$ so that there exists a chain $w_1\prec_\fw w_2\hdots\prec_\fw w_n$ and write $n={\rm hgt}(\fw)$.
The following is then immediate from Lemma \ref{LemmNotTall}.

\begin{corollary}\label{CorNotTall}
If $\Sigma\Subset \mathcal L$ and $\fw$ is an irreducible $\Sigma$-moment, then ${\rm hgt}(\fw)\leq \#\Sigma+1$.
\end{corollary}

Next we will give a bound on the size of irreducible moments. Our bound will be superexponential; recall that the superexponential $2^m_n$ is defined by recursion on $n$ by $2^m_0=m$ and $2^m_{n+1}=2^{2^m_{n}}$. We begin with a useful inequality:

\begin{lem}
For all $m,n,k\geq 1$, $2^m\cdot 2_{k}^{(n-1)m}\leq 2_{k}^{nm}.$
\end{lem}

\proof
Proceed by induction on $k$. If $k=1$, then
\[2^m\cdot 2_{k}^{(n-1)m}=2^m\cdot 2^{(n-1)m}=2_{k}^{n m}.\]
If $k>1$, then note that $1 \leq  2_{k-1}^{(n-1)m}$, so that
\[m+ 2_{k-1}^{(n-1)m} \leq (m+1) 2_{k-1}^{(n-1)m} \leq 2^m\cdot 2_{k-1}^{(n-1)m} \stackrel{\text{{\sc IH}}} \leq 2_{k-1}^{nm} .\]
Then,
\[2^m\cdot 2_{k}^{(n-1)m}=2^m \cdot 2^{ 2_{k-1}^{(n-1)m}}=2^{m+ 2_{k-1}^{(n-1)m}} \leq 2^{2_{k-1}^{nm}} = 2_{k}^{nm},\]
as needed.
\endproof

\begin{lem}\label{LemIrrBound}
If $\Sigma\Subset\landi$ and $s=\#\Sigma>0$, then 
$\# I_\Sigma \leq 2^{s^2+s}_{s+1}.$
\end{lem}

\proof
We prove, by induction on $n$, that there are at most $2^{n s}_{n}$ irreducible moments with height $n$. Note that there are no moments of height $0$, and $0=2^{0s}_0$. So, we may assume $n\geq 1$. Any irreducible of height at most $n$ is of the form $\add \Phi A$, where $\Phi\subseteq \Sigma$ and $A$ is a set of distinct irreducibles of height less than $n$. There are at most $2^{s}$ choices for $\Phi$. By induction hypothesis, there are at most $2^{(n-1)s}_{n-1}$ possible elements for $A$, and thus $2^{2^{(n-1)s}_{n-1}}=2^{(n-1)s}_{n}$ choices for $A$. Hence there are at most
$2^{s}\cdot 2^{(n-1)s}_{n}\leq 2^{n s}_{n}$ such $\fw$. The lemma then follows by choosing $n=s+1$, the greatest value that ${\rm hgt}(\fw)$ could take.
\endproof

Now that we know that $\irr\Sigma$ is finite, it would be convenient if, whenever $\chi$ is a simulation and $\fw$ is any $\Sigma$-moment such that $\fw\mathrel \chi x$, we could replace $\fw$ by some irreducible $\fw'\redu\fw$ and still have $\fw'\mathrel \chi x$. The following operations on simulations will help us achieve this.

\begin{defn}
Let $\Sigma\Subset \Delta \subseteq \landi$ both closed under subformulas, $\mathfrak X$ be a $\Delta$-labelled space, and $\chi\subseteq M_\Sigma\times |\mathfrak X|$.
\begin{enumerate}

\item Define the {\em reductive closure} of $\chi$ by $\check \chi=\chi{ \redu }$. If $\check \chi=\chi$, we say $\chi$ is {\em reductive.}

\item Define the {\em irreducible part} of $\chi$ by $\chi_0=\chi\upharpoonright I_\Sigma$.

\end{enumerate}

\end{defn}

In other words, $\fw \mathrel{\check \chi} x$ means that there is $\fv$ such that $\fw \redu \fv \mathrel \chi x$. Let us see that these operations indeed produce new simulations.

\begin{lem}\label{LemRedSims}
Suppose that $\Sigma\Subset\landi$, $\mathfrak X$ is a labelled system, and $\chi\subseteq M_\Sigma\times |\mathfrak X|$ is a simulation. Then,

\begin{enumerate}

\item $\check \chi$ and $\chi_0$ are also simulations.

\item If $\chi$ is reductive, then
\begin{enumerate}

\item\label{ItRedSimsA} $\chi_0(I_\Sigma)=\chi(M_\Sigma)$, and

\item\label{ItRedSimsB} if $\mathfrak w\mathrel\chi x \mathrel S_\mathfrak X y$, $\fv\mathrel \chi y$, and $\mathfrak w\circrel\mathfrak v$, there is $\mathfrak v'$ such that $\mathfrak v'\mathrel \chi_0 y$ and $\mathfrak w\circrel \mathfrak v'$.

\end{enumerate}

\end{enumerate}

\end{lem}

\proof
That $\check\chi$ is a simulation is immediate from Lemma \ref{LemmPropSim}.\ref{ItPropSimComp}, and it is easily seen that $\chi_0$ is a simulation using Lemma \ref{LemmPropSim}.\ref{ItPropSinSub}.

For claim \ref{ItRedSimsA}, it is obvious that $\chi_0(I_\Sigma)\subseteq \chi(M_\Sigma)$. For the other inclusion, assume that $\chi$ is reductive. Consider $x\in \chi(M_\Sigma)$, so that there is $\fw\in M_\Sigma$ with $\fw\mathrel \chi x$. By Lemma \ref{LemRedIsSim}, $\redu_0$ is surjective, so we can pick $\fw'$ with $\fw'\redu_0\fw$; since $\chi$ is reductive, it follows that $\fw'\mathrel\chi x$, and hence $x\in \chi_0(I_\Sigma)$, as needed.

For Claim \ref{ItRedSimsB}, choose $\fv'\redu_0\fv$, so that, as above, $\fv'\mathrel \chi_0 y$. Meanwhile, since $\fw\redu \fw$ and $\fw \mathrel S_\Sigma \fv$, we can use Lemma \ref{LemPropRed}.\ref{ItLemPropRedNext} to conclude that $\fw \mathrel S_\Sigma \fv'$.
\endproof

\begin{figure}

\begin{center}

\begin{tikzpicture}[scale = 0.7]

\draw[thick] (0,0) circle (.35);

\draw (.03,-.04) node {$x$};

\draw[thick,->] (.5,0) -- (3.5,0);

\draw (2,-.5) node {$S_\mathfrak X$};

\draw[thick] (4,0) circle (.35);

\draw (4+.03,-.05) node {$y$};

\draw[thick] (0,4) circle (.35);

\draw (.04,4-.04) node {$\fw$};

\draw[thick,->] (.5,4) -- (3.5,4);

\draw (2,4.5) node {$S_\Sigma$};

\draw[thick,->] (0,3.5) -- (0,.5);

\draw (-.5,2) node {$\chi$};

\draw[thick] (4,4) circle (.35);

\draw (4+.05,4-.04) node {$\fv$};

\draw[thick,->] (4,3.5) -- (4,.5);

\draw (1.2,3.3) node {$S_\Sigma$};

\draw[thick,->,dashed] (.4,4-.4) -- (2-.35,2+.35);

\draw (4.5,2) node {$\chi$};

\draw[thick,dashed] (2,2) circle (.35);

\draw (2+.03,2-.04) node {$\fv'$};

\draw (2.8,.7) node {$\chi_0$};

\draw[thick,<-,dashed] (4-.4,.4) -- (2+.35,2-.35);

\end{tikzpicture}

\end{center}

\caption{A diagram illustrating Lemma \ref{LemRedSims}.\ref{ItRedSimsB}.}

\end{figure}
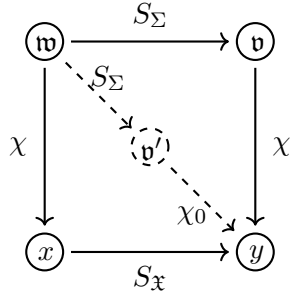

\section{Simulating Dynamical Systems}\label{SecSDS}

Our work with irreducibles shows that if we can construct a surjective, dynamic simulation on $\moments\Sigma$, then we immediately get a surjective, dynamic simulation on $\irr\Sigma$, which we can then use to construct finite quasimodels. In this section, we will show how such simulations can be found; specifically, we will show that maximal simulations have all required properties.

\begin{lem}
Given any labelled space
$\mathfrak{X}$ and $\Sigma\Subset\landi$,
there exists a greatest simulation
$\chi^\ast\subseteq M_\Sigma\times |\mathfrak X|.$ Moreover, $\chi^\ast$ is reductive.
\end{lem}
\proof
By Lemma \ref{LemmPropSim}.\ref{ItPropSimUn},
\[\chi^\ast=\bigcup\{\chi\subseteq M_ \Sigma \times |\mathfrak X|: \chi\text{ is a simulation} \}\]
is a simulation, and clearly this is the greatest simulation. By Lemma \ref{LemRedSims}, $\check\chi^\ast$ is also a simulation, hence $\check\chi^\ast\subseteq \chi^\ast$, and thus $\chi^\ast$ is reductive.
\endproof

If $\mathfrak X$ is a labelled system then $\chi^\ast$ gives us a surjective, dynamic simulation. The following lemma will be essential in proving this.
\begin{lem}\label{simex}
Let $\Sigma\Subset\landi$,
$\mathfrak X$
be a labelled space, and
$\chi\subseteq M_\Sigma\times |\mathfrak X|$
be a simulation.

Suppose that $\Phi\subseteq \Sigma$ and $A\subseteq M_\Sigma$ are a $\Sigma$-kit (as in Definition \ref{cohe}), and that there is $x_\ast\in |\mathfrak X|$ such that
\begin{enumerate*}

\item $\ell_\mathfrak X(x_\ast) \cap \Sigma =\Phi$, and

\item
for all $\fw\in A$,
$x_\ast \in \overline{\chi(\fw)}$.
\end{enumerate*}
Then, $\zeta=\chi\cup\cbra\left (\add \Phi A,x_\ast \right )\cket$ is also a simulation.
\end{lem}
\proof
Note, first, that $\zeta$ is label-preserving; let us prove that it is also continuous.
Consider an arbitrary open set $U\subseteq |\mathfrak X|$. We must show that $\zeta^{\,-1}(U)$ is open. So, let $\fw\in \zeta^{\,-1}(U)$, and assume that $\fv\acc\fw$; we must show that $\fv\in \zeta^{\,-1}(U)$. Obviously the latter holds if $\fw=\fv$, so we assume $\fv\prec \fw$. Pick $x\in U$ such that $\fw\mathrel\zeta x$. If $x\not=x_\ast$ or $\fw\not=\add\Phi A$, then we already have $\fw\mathrel \chi x$, and hence $\fv\in \chi^{\,-1}(U)\subseteq \zeta^{\,-1}(U)$. So, assume that $x=x_\ast$ and $\fw=\add\Phi A$. Then, there is $\fv'\in A$ such that $\fv\acc \fv'$. By assumption, $x_\ast \in \overline{\chi(\fv')}$; hence, since $U$ is open, there is $y\in U$ such that $\fv' \mathrel \chi y$. But $\chi$ is continuous and $U$ is a neighborhood of $y$, so there is $z\in U$ such that $\fv \mathrel \chi z$, hence $\fv\mathrel \zeta z$ as needed.
\endproof

\begin{prop}\label{total}
Let $\Sigma\Subset \Delta \subseteq \landi$ both be closed under subformulas. Then, for any $\Delta$-labelled space $\mathfrak{X}$, $\chi^\ast\subseteq M_\Sigma\times |\mathfrak X|$ is surjective.
\end{prop}
\proof
Since we know that $\chi^\ast$ is reductive, it suffices to show that a reductive simulation that is not surjective is not maximal. Thus, let
$\chi\subseteq M_\Sigma\times |\mathfrak X|$
be any reductive simulation. Say a point $x\in X$ is {\em bad} if it does not lie in the range of $\chi$. We claim that $\chi$ is not maximal if there are bad points.

If there {\em are} bad points, let $x_\ast\in |\mathfrak X|$ be a bad point such that $\Phi := \ell(x_\ast) \cap \Sigma $ is $\subseteq$-maximal among all bad points. Let $U_\ast$ be a neighborhood of $x_\ast$ minimizing $\#\chi^{-1}_0(U_\ast)$ (which is finite by Lemma \ref{LemIrrBound}), and such that $U_\ast\subseteq \lb\varphi\rb_\mathfrak X$ for all $\varphi\in \Phi $. We claim that for each defect $\delta= (\varphi\imp\psi) \in\defect\Phi$, there is $\mathfrak u_\delta\in \chi^{-1}_0(U_\ast)$ revoking $\delta$. To see this, note that by the semantics of $\imp$, there is $x_\delta\in U_\ast$ such that $x_\delta\in \lb\varphi\rb_\mathfrak X\setminus \lb\psi\rb_\mathfrak X$. But, since $x_\delta\in U_\ast$, it follows that $\ell(x_\delta) \cap \Sigma \supseteq \Phi $, and since $\varphi \in \ell(x_\delta)\setminus \ell(x_\ast)$, the inclusion is strict. By maximality of $\Phi$ it follows that $x_\delta$ is not bad, hence $x_\delta\in \chi(\mathfrak v)$ for some $\mathfrak v\in M_\Sigma$. By Lemma \ref{LemExistIrr}, there is an irreducible $\mathfrak u_\delta\redu\mathfrak v$, and since $\chi$ is reductive, $\mathfrak u_\delta\mathrel\chi x_\delta$, so that $\mathfrak u_\delta\in \chi^{-1}_0(U_\ast)$. By the minimality of $\#\chi^{-1}_0(U_\ast)$ we see that if $U$ is any neighborhood of $x$ then $\mathfrak u_\delta \in \chi^{-1}_0(U \cap U_\ast) \subseteq \chi^{-1}_0(U )$, so that $x_\ast \in \overline{\chi(\mathfrak u_\delta)}$.

Let $A$ be the set of all $\mathfrak u_\delta$ such that $\delta \in \partial \Phi$. By Lemma \ref{add},
$\fw_\ast=\add{\Phi }{A}$ is a $\Sigma$-moment, and we can set
\[\zeta=\chi\cup \big\{ \big (\textstyle\add{\Phi }{A} , x_\ast \big)\big \}.\]
By Lemma \ref{simex}, $\zeta$ is a simulation. Since $x_\ast$ was bad, $\chi\subsetneq \zeta$, as desired.
\endproof

\begin{prop}\label{temptotal}
Given $\Sigma\Subset \Delta \subseteq \landi$ both closed under subformulas and a $\Delta$-labelled system $\mathfrak{X}$, $\chi^\ast \subseteq M_\Sigma \times |\mathfrak X|$ is a dynamic simulation.
\end{prop}

\proof
The proof follows much the same structure as that of Proposition \ref{total}. Let $\chi$ be a reductive simulation, which in view of Proposition \ref{total}, we may assume to be surjective. Suppose further that $\chi$ is not dynamic; we will show that it cannot be maximal.

Say that $(x,y) \in S_\mathfrak X$ {\em fails} for $\mathfrak w$ if $\mathfrak w\mathrel \chi x$ but there is no $\mathfrak v\in \chi^{-1} (y)$ such that $\mathfrak w \circrel \mathfrak v$. We will show that $\chi$ is not maximal if any pair fails for any moment. If this were the case, pick $\mathfrak w_\ast$ of minimal height such that some pair $(x_\ast, y_\ast ) \in S_\mathfrak X$ fails for it. Let $V_\ast$ be a neighborhood of $ y_\ast $ such that $\# \chi^{-1}_0(V_\ast)$ is minimal and $V_\ast\subseteq \lb\varphi\rb_\mathfrak X$ whenever $\varphi\in \Psi := \ell(y_\ast) \cap \Sigma$. 

As before, for each defect $\delta$ of $\ell(y_\ast )$, choose $\mathfrak u_\delta\in \chi^{-1}_0(V_\ast)$ revoking $\delta$. Next, for each $\mathfrak v$ such that $\fv\prec \fw$, we claim there is $\mathfrak v'\in \chi^{-1}_0(V_\ast)$ such that $\mathfrak v\circrel\mathfrak v'$. To see this, using the continuity of $S_\mathfrak X$, let $U$ be a neighborhood of $x_\ast$ such that $U\subseteq S^{-1}_\mathfrak X (V_\ast)$. Since $\chi$ is continuous, there is $x_\mathfrak v\in U$ such that $\mathfrak v\mathrel \chi x_\mathfrak v$, and by our choice of $U$ there is $y_\mathfrak v \in V_\ast$ such that $x_\mathfrak v \mathrel S_\mathfrak X y_\mathfrak v$. Since ${\rm hgt}(\mathfrak v)< {\rm hgt}(\fw)$, by the induction hypothesis, there is $\mathfrak v''$ such that $\mathfrak v\circrel\mathfrak v''  \mathrel \chi y_\fv$. But $\chi$ is reductive, so that by Lemma \ref{LemRedSims}.\ref{ItRedSimsB}, there is $\fv'$ such that $\fv'\mathrel \chi_0 y_\mathfrak v \in V_\ast$ and $\fv\circrel\fv'$, as needed.

Let
\[B= \big \{\fu_\delta: \delta\in\defect{\ell(y_\ast)}  \big \}\cup \{\mathfrak v':\mathfrak v\prec\mathfrak w\}.\]
By the minimality of $\# \chi^{-1}_0(V_\ast)$ we see that $y_\ast \in \overline {\chi(\mathfrak b)}$ for all $\mathfrak b \in B$.
By Lemma \ref{add}, $\mathfrak v_\ast=\add{\Psi }B$ is a $\Sigma$-moment and $\mathfrak w_\ast\circrel\mathfrak v_\ast$ by Lemma \ref{add}.\ref{ItAddTwo}.

We can then set
\[\zeta=\chi\cup \cbra (\mathfrak v_\ast,y_\ast)\cket.\]
As before, $\zeta$ is a simulation which properly contains $\chi$. Therefore, $\chi$ is not maximal.
\endproof

We are now ready to prove that $\icltl$ has the effective finite quasimodel property and hence is decidable.

\begin{defn}
If $\Sigma \Subset \Delta \subseteq \landi$ are both closed under subformulas and $\mathfrak X$ is a $\Delta$-labelled space, we define
$\mathfrak X/\Sigma=\mathbb I_\Sigma\upharpoonright \mathrm{dom}(\chi^\ast).$
\end{defn}

\begin{prop}\label{PropQuasiComplete}
Let $\Sigma \Subset \Delta \subseteq \landi$ be closed under subformulas and $\mathfrak X$ be a well $\Delta$-labelled system.
Then, if $\mathfrak X$ falsifies $\varphi\in \Sigma$, it follows that $\mathfrak X/\Sigma$ is a $\Sigma$-quasimodel falsifying $\varphi$.
\end{prop}

\proof
By Propositions \ref{total} and \ref{temptotal}, $\chi^\ast$ is a surjective dynamic simulation, so by Lemmas \ref{LemTotSim} and \ref{isquasi}, ${\mathfrak X}/\Sigma$ is a $\Sigma$-quasimodel. Now, fix $\varphi \in \Sigma$ and $x \in |\mathfrak X| $. By Proposition \ref{total}, $\chi^\ast_0$ is surjective, so if there is $x \in |\mathfrak X| \setminus \val{\varphi}_\mathfrak X$ there exists $\fw \in I_\Sigma$ such that $x \in\chi^\ast(\fw )$ and hence $\varphi \not \in \ell(\fw )$.
\endproof

This gives us a finite quasimodel property for $\icltl$.

\begin{thm}
A formula $\varphi \in \landi$ is falsifiable on a dynamical system if and only if it is falsifiable on a quasimodel of size at most $2^{(\lgt\varphi+1)\lgt\varphi}_{\lgt\varphi}$.
\end{thm}

\proof
Let $\varphi \in \landi$ and $\Sigma={\rm sub}(\varphi)$. By Proposition \ref{second}, if $\varphi$ is falsifiable on a $\Sigma$-quasimodel, then it is falsifiable on a dynamic topological model. Conversely, if $\varphi$ is falsifiable on a dynamic topological model $\mathfrak X$, then by regarding $\mathfrak X$ as a well $\landi$-labelled system we see using Proposition \ref{PropQuasiComplete} that $\varphi$ is falsifiable on a $\Sigma$-quasimodel of the form $\mathfrak Q=\mathfrak X/\Sigma$, which is a substructure of $\mathbb I_\Sigma$. Since $\mathfrak Q$ has at most $2^{(\lgt\varphi+1)\lgt\varphi}_{\lgt\varphi}$ worlds, our result follows.
\endproof

From this our main result is immediate.

\begin{thm}\label{TheoITLDec}
$\icltl$ is decidable.
\end{thm}

\section{Special Classes of Systems}\label{SecSpec}

Now that we have seen that $\icltl$ is decidable, let us turn our attention to its capacity to reason about recurrence phenomena. The operator $\ps$ allows us to express properties of the asymptotic behavior of dynamical topological systems; this behavior is often more interesting when these systems have additional structural properties. We will discuss two special classes which exhibit non-trivial asymptotic behavior, and show that $\icltl$ can describe this behavior in both cases. We begin by discussing minimal systems.

\subsection{Minimal systems}\label{SecMinimal}

Minimal systems, introduced in \cite{Bir}, are dynamical systems that exhibit characteristic asymptotic behavior which can be captured in the language of ${\sf DTL}$. We will recall a few basic properties for minimal dynamical systems; a deeper treatment can be found in a text such as \cite{akin}. 

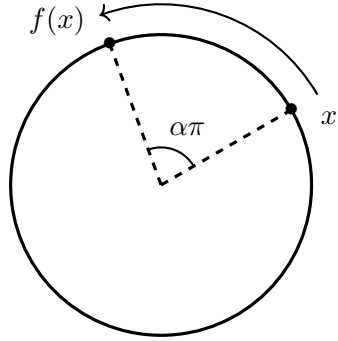
\begin{figure}

\def \cosine {0.86602540378443864676372317075294}

\def \sine {0.5}

\def \costwo {-0.34202014332566873304409961468226}

\def \sintwo {0.93969262078590838405410927732473}

\begin{center}

\begin{tikzpicture}

\draw[thick] (.5*\cosine,.5*\sine) arc (30:110:.5);

\draw[very thick] (0,0) circle (2);

\draw[very thick,dashed] (0,0) -- (2*\cosine,2*\sine);

\draw[very thick,dashed] (0,0) -- (2*\costwo,2*\sintwo);

\draw (2*\cosine,2*\sine) node {$\bullet$};

\draw (2*\costwo-.7,2*\sintwo+.3) node {$f(x)$};

\draw(2*\costwo,2*\sintwo) node {$\bullet$};

\draw(2*\cosine+.5,2*\sine-.1) node {$x$};

\draw (.35,.75) node {$\alpha\pi$};

\draw[thick,->] (1.2*2*\cosine,1.2*2*\sine) arc (30:110:1.2*2);

\end{tikzpicture}

\end{center}

\caption{If $\alpha$ is irrational, then a rotation by $\alpha \pi$ gives rise to a minimal system on the unit circle.}
\label{FigCirc}

\end{figure}

If $\mathfrak X$ is a dynamical system and $Y\subseteq|\mathfrak X|$ is such that $ f_\mathfrak X (Y) \subseteq Y$, then $\mathfrak X\upharpoonright Y$ denotes the dynamical system
\[\<Y,\mathcal T_\mathfrak X \upharpoonright Y,f_\mathfrak X\upharpoonright Y\>,\]
where $\mathcal T_\mathfrak X \upharpoonright Y$ is the subspace topology on $Y$ induced by $\mathcal T_\mathfrak X$; to be precise,
\[(\mathcal T_\mathfrak X \upharpoonright Y) = \{U\cap Y : U \in \mathcal T_\mathfrak X\}.\]

\begin{defn}
Let $\mathfrak X$ and $\mathfrak Y$ be dynamical systems. We will say that $\mathfrak Y$ is a {\em subsystem} of $\mathfrak X$, written $\mathfrak Y\leq\mathfrak X$, if there exists a non-empty closed set $Y\subseteq|\mathfrak X|$ which is closed under $f_{\mathfrak X}$ and such that $\mathfrak Y=\mathfrak X\upharpoonright Y$.

The dynamical system $\mathfrak X$ is {\em minimal} if it is minimal under $\leq$; that is, if it contains no proper subsystems. 
\end{defn}

While it is not true in general that a dynamical system $\mathfrak X$ has a minimal subsystem, the following theorem of Birkhoff~\cite{Bir} shows that this is the case if we assume $\mathfrak X$ to be compact.

\begin{thm}[Birkhoff]
If $\mathfrak X$ is a dynamical system based on a compact space, then there exists a minimal $\mathfrak Y\leq\mathfrak X$.
\end{thm}

The result is proven using Zorn's lemma. An alternative, and sometimes more useful, characterization of minimal systems is given by the following:

\begin{prop}
A dynamical system $\mathfrak X$ is minimal if and only if the orbit of every point is dense in $|\mathfrak X|$.
\end{prop}

\proof
Suppose that $\mathfrak X$ is minimal and let $x\in|\mathfrak X|$. Let $Y$ be the the closure of the orbit of $x$:
\[Y=\overline{\cbra f^n_{\mathfrak X}(x):n<\omega\cket}.\]
It is not hard to check that $Y$ is a closed subset of $\mathfrak X$ that is closed under $f_{\mathfrak X}$; since $\mathfrak X$ is minimal, it follows that $Y=|\mathfrak X|$, which means that the orbit of $x$ is dense in $|\mathfrak X|$.

Now, suppose that $\mathfrak X$ is not minimal, and let $Z\subsetneq |\mathfrak X|$ be non-empty, closed, and $f_{\mathfrak X}$-closed. Note that $Z$ itself cannot be dense since it is already closed. Then, it is clear that for any $x\in Z$, the orbit of $x$ lies in $Z$, and since $Z$ is closed, the closure of the orbit of $x$ is also a subset of $Z$. Thus the orbit of $x$ is not dense in $|\mathfrak X|$.
\endproof

An example of a minimal system is given by a rotation of the unit circle $S^1$ by an angle $\alpha \pi$, where $\alpha$ is irrational. It is well-known that under such a rotation, the orbit of every point of $S^1$ is dense, and hence this gives us a minimal system (see Figure \ref{FigCirc}). Note that minimal systems do not have to be connected or surjective; Figure \ref{FigMinNonInv} gives an example of a minimal system that is neither, and compact, Hausorff examples can be found in \cite{Balibrea,Bruin}.

The dynamic topological logic of minimal systems is well-understood \cite{me:minimal}:

\begin{thm}\label{TheoDTLMinExpress}
Let $\mathfrak X$ be a dynamical system. Then, $\mathfrak X$ is minimal if and only if
\[\mathfrak X\models\exists\nb p\to\forall\ps p.\]
\end{thm}

\begin{thm}\label{TheoDTLMin}
Let $\DTLm$ denote the set of valid formulas of $\lanclass$ over the class of minimal systems. Then, $\DTLm$ is decidable.
\end{thm}

As we will see, analogous results hold for $\lanfull$. First, observe that decidability is readily inherited from the classical setting:

\begin{thm}\label{TheoITLMin}
Let $\ITLm$ denote the set of valid formulas of $\lanfull$ over the class of minimal systems. Then, $\ITLm$ is decidable.
\end{thm}

\proof
We know that $\varphi\in \lanfull$ is valid over the class of minimal systems if and only if $\varphi^\blacksquare$ is. Thus, the validity problem of $\ITLm$ can be reduced to that of $\DTLm$; since the latter is decidable, so is the former.
\endproof

\begin{exam}\label{ExMinCounter}
Figure \ref{FigMinNonInv} gives an example of an Alexandroff dynamical system $\mathfrak m$, where $f_\mathfrak m$ is not open. Note that the set $\{w,z\}$ is dense, and the orbit of every point contains these two points, so the system is minimal. If we let $\val p=\{y,z\}$ and $\val q=\{z\}$, then it is readily verified that $v\not\in \val{(\nx p \imp \nx q) \imp \nx (p\imp q)}$, hence $(\nx p \imp \nx q) \imp \nx (p\imp q) \not \in \ITLm$.
 
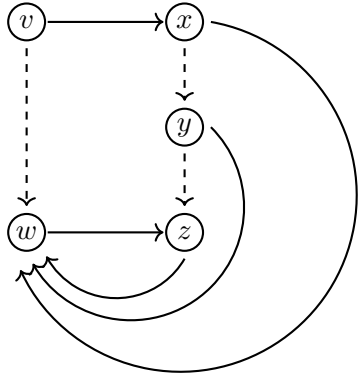
\begin{figure}
\begin{center}

\begin{tikzpicture}[scale=0.7]

\draw[thick] (0,0) circle (.35);

\draw (0,0) node {$w$};

\draw[thick,->,dashed] (0,3.5) -- (0,0.5);

\draw[thick] (0,4) circle (.35);

\draw (0,4) node {$v$};

\draw[thick] (3,0) circle (.35);

\draw (3,0) node {$z$};

\draw[thick,->,dashed] (3,1.5) -- (3,.5);

\draw[thick] (3,2) circle (.35);

\draw (3,2) node {$y$};

\draw[thick,->,dashed] (3,3.5) -- (3,2.5);

\draw[thick] (3,4) circle (.35);

\draw (3,4) node {$x$};

\draw[thick,-> ] (.4,0) -- (2.6,0);

\draw[thick,-> ] (.4,4) -- (2.6,4);

\draw[thick,->] (3,-0.5) arc (-30:-152:1.5) ;

\draw[thick,->] (3.5,2) arc (45:-150:2.15) ;

\draw[thick,->] (3.5,4) arc (80:-155:3.35) ;

\end{tikzpicture}

\end{center}
\caption{An Alexandroff minimal system $\mathfrak m$ with the down-set topology. Dashed arrows indicate $\peq$, solid arrows $f_\mathfrak m$.}\label{FigMinNonInv}
\end{figure}

\end{exam}

\begin{remark} The set of valid formulas of $\DTLm$ is not primitive recursively decidable, hence the above argument does not settle whether $\ITLm$ is. One could combine the proof of Theorem \ref{TheoITLDec} with the proof in \cite{me:minimal} of Theorem \ref{TheoDTLMin} to obtain a primitive recursive decision procedure; nevertheless, the procedure thus obtained would most likely remain non-elementary without some non-trivial modifications.
\end{remark}

As for characterizing minimality in $\lanfull$ (indeed, in $\landi$), we may actually do this more succinctly than in $\lanclass$:

\begin{thm}\label{TheoITLMinExpress}
Let $\mathfrak X$ be a dynamical system. Then, $\mathfrak X$ is minimal if and only if
\begin{equation}\label{EqMinFor}
\mathfrak X\models\exists  p\imp \ps p.
\end{equation}
\end{thm}

\proof
First assume that $\mathfrak X$ is any minimal dynamical system and $\val \cdot$ is any valuation on $\mathfrak X$. Suppose that $x\in \val{\exists p}$. Then, $\val p $ is non-empty and open; since the orbit of $x$ is dense, it follows that $f^n_\mathfrak X(x) \in \val p $ for some $n$, that is, $x\in \val{\ps p}$.

Conversely, if $\mathfrak X$ is {\em not} minimal, let $Y$ be a non-empty, proper, closed, $f_\mathfrak X$-closed subset of $|\mathfrak X|$. Then, if $\val\cdot$ is any valuation satisfying $\val p=|\mathfrak X|\setminus Y$ and $x\in Y$, it is not hard to check that $x\not\in \val{\exists p \imp \ps p}$.
\endproof

As a corollary, we obtain the following:

\begin{corollary}
Over the class of minimal systems, $\ps$ is definable by $\ps \varphi \iiff \exists \varphi$.
\end{corollary}

In view of Lemma \ref{LemmExDefin}, we may also define $\ps \varphi \iiff \ineg \forall \ineg \varphi$, and thus the language $\lanfull\upharpoonright \{\nx,\forall\}$ is equally expressive as $\landi$. This makes $\lanfull\upharpoonright \{ \nx, \forall\}$ (i.e., the language with modalities $\nx,\forall$), or even $\lanfull\upharpoonright \{ \nx, \exists \}$, be particularly attractive for reasoning about minimal dynamical systems.

\subsection{Poncar\'e recurrence}\label{SecProb}

Finally, we turn our attention to probability-preserving systems. The $\sf DTL$ of such systems is not as well-understood as that of minimal systems; nevertheless, the Poincar\'e recurrence theorem, one of the motivations for the study of $\sf DTL$, holds in this class of systems.

\begin{defn}
A dynamical system $\mathfrak X$ is {\em probabilty preserving} if there exists a probability measure $\pi$ on $|\mathfrak X|$ such that
\begin{enumerate}
\item if $U\subseteq|\mathfrak X|$ is open and non-empty, then $U$ is measurable and $\pi(U)\not=0$;
\item for every measurable set $Y\subseteq|\mathfrak X|$, we have $\pi(Y)=\pi(f^{-1}_\mathfrak X(Y))$.
\end{enumerate}
\end{defn}

\begin{figure}

\def \cosine {0.86602540378443864676372317075294}

\def \sine {0.5}

\def \costwo {-0.34202014332566873304409961468226}

\def \sintwo {0.93969262078590838405410927732473}

\begin{center}

\begin{tikzpicture}[scale=0.8]

\draw[very thick] (0,0) circle (3);

\draw[thick] (.5*\cosine,.5*\sine) arc (30:110:.5);

\draw[very thick,dashed] (0,0) circle (2);

\draw[very thick,dashed] (0,0) -- (3*\cosine,3*\sine);

\draw[very thick,dashed] (0,0) -- (3*\costwo,3*\sintwo);

\draw (2*\cosine,2*\sine) node {$\bullet$};

\draw (2*\costwo-.7,2*\sintwo+.1) node {$f(x)$};

\draw(2*\costwo,2*\sintwo) node {$\bullet$};

\draw(2*\cosine+.5,2*\sine-.1) node {$x$};

\draw (.35,.75) node {$\vartheta$};

\draw[very thick,->] (2.2*\cosine,2.2*\sine) arc (30:109:2.2);

\end{tikzpicture}

\end{center}

\caption{A rotation $f$ of a disk by an angle $\vartheta$ gives rise to a probability-preserving system, where the probability of a set is defined to be proportional to its area. Note that $\vartheta$ no longer needs to be an irrational multiple of $\pi$, unlike in Figure \ref{FigCirc}. Moreover, such a system is not minimal, as (for example) the dashed circle is a closed, $f$-closed subsystem.}
\label{FigPoinc}

\end{figure}
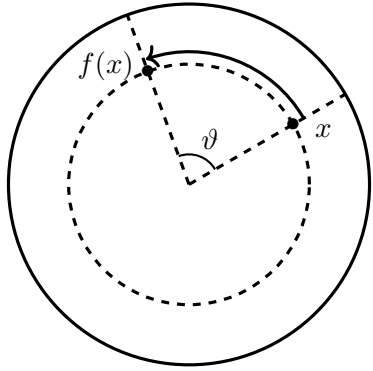

Probability-preserving systems have the following property, discovered by Poincar\'e \cite{Poincare}; a more modern presentation can be found in \cite{akin}.

\begin{thm}\label{TheoPoincare}
If $\mathfrak X$ is a probability-preserving system and $U \subseteq|\mathfrak X|$ is open, then $U$ contains a recurrent point, i.e.~there is $x\in U$ such that $f^n_\mathfrak X (x)\in U$ for some $n>0$.
\end{thm}

As observed by Kremer and Mints in \cite{kmints}, Poincar\'e recurrence can readily be expressed in the language of $\sf DTL$.

\begin{thm}\label{TheoDTLPoincare}
The $\lanclass$ formula
\begin{equation}
\nb p \to \pb \nx \ps p \label{EqPoincare}
\end{equation}
is valid over the class of probability-preserving systems, but not valid in general.
\end{thm}

Poincar\'e recurrence is also expressible in $\landi$. Below, say that a dynamical system $\mathfrak X$ is {\em Poincar\'e recurrent} if every non-empty open set $U\subseteq |\mathfrak X|$ contains a recurrent point.

\begin{thm}\label{TheoITLPoincare}
A dynamical system $\mathfrak X$ is Poincar\'e recurrent if and only if
\begin{equation}\label{EqPoincInt}
\mathfrak X\models  p \imp \ineg \ineg   \nx \ps p .
\end{equation}
\end{thm}

\proof
First assume that $\mathfrak X$ is Poincar\'e recurrent, and let $\val\cdot$ be any valuation on $\mathfrak X$. Let $U\subseteq \val p$ be open and non-empty. Then, $U$ contains a recurrent point; that is, there are $x\in U$ and $n >0 $ such that $f^n_\mathfrak X(x) \in U$. It follows that $x \in \val{ \nx\ps p}$, and since $U\subseteq \val p$ was arbitrary, using Lemma \ref{LemmDoubleNeg} we conclude that $\mathfrak X \models p \imp \ineg \ineg  \nx \ps p $.

Conversely, if $\mathfrak X$ is not Poincar\'e recurrent, let $U \subseteq |\mathfrak X|$ be a non-empty open set with no recurrent point. Set $\val p = U$ and choose $x\in U$; then, it is easy to check that $x \not\in \val {\nx \ps p}$. Since $x\in U$ is arbitrary we conclude that $U\cap \val {\nx \ps p} =\varnothing$, hence $\val {\nx \ps p} =\varnothing$ is not dense in $\val p$, and $\mathfrak X\not\models p\imp \ineg\ineg \nx\ps p$.
\endproof

\begin{remark}
Despite being able to capture minimality and Poincar\'e recurrence, $\lanfull$ has less expressive power than $\lanclass$, since the former can only reason about open sets. Thus, we should expect the intuitionistic language to have some limitations. One example of this might be given by {\em strong} Poincar\'e recurrence.

Say that a dynamical system $\mathfrak X$ is {\em strongly Poincar\'e recurrent} if every non-empty open $U\subseteq |\mathfrak X|$ has an {\em infinitely recurrent} point; that is, there is $x \in U$ such that $f^n_\mathfrak X (x) \in U$ for infinitely many $n$. Then, Theorem \ref{TheoPoincare} can be strengthened to conclude that any probability-preserving system is strongly Poincar\'e recurrent. 

One can then reason as in the proof of Theorem \ref{TheoITLPoincare} to conclude that $\mathfrak X \models \nb p \to \pb \nec \ps p$ if and only $\mathfrak X$ it is strongly Poincar\'e recurrent. However, given that the set of infinitely recurrent points need not be open, it is unlikely that $\lanfull$ can characterize strongly Poincar\'e recurrent systems in this fashion.
\end{remark}

\section{Concluding Remarks}\label{SecConc}

\begin{table}

\begin{center}

\begin{tabular}{| p{4cm} | p{1.5cm} | p{4cm} | p{1.7cm} | p{1.7cm} |}

\hline

\multicolumn{1}{|c|}{Class} & \multicolumn{1}{c|}{Notation} &\multicolumn{1}{c|}{$\phantom{\displaystyle{a\choose b}}\lanclass\phantom{\displaystyle{a\choose b}}$} & \multicolumn{1}{c|}{$\landi$} & \multicolumn{1}{c|}{$\lanfull$}\\

\hline

All dynamical systems & \multicolumn{1}{c|}{$\sf c$} & Undecidable \cite{konev} but c.e. \cite{me2} & Decidable & Unknown but c.e.\\

\hline

Dynamical systems with a homeomorphism & \multicolumn{1}{c|}{$\sf h$} & Non-c.e. \cite{wolter} & Unknown & Unknown\\

\hline

Expanding frames & \multicolumn{1}{c|}{$\sf e$} & Undecidable \cite{konev} & Decidable  &  Decidable \cite{Boudou2017} \\

\hline

Persistent frames & \multicolumn{1}{c|}{$\sf p$} & Non-c.e. \cite{wolter} & Unknown  & Unknown \phantom{but c.e.} \\

\hline

Minimal systems & \multicolumn{1}{c|}{$\sf m$} & Decidable, but not primitive recursive \cite{me:minimal} & Decidable  & Decidable \\

\hline

Poincar\'e recurrent systems & \multicolumn{1}{c|}{$\sf r$} & Unknown & Unknown  & Unknown \\

\hline

\end{tabular}

\end{center}

\caption{This table indicates whether the set of formulas of a given language is decidable over different classes of dynamical systems. Note that the notations $\sf h$, $\sf r$ have not been introduced previously in the text.}\label{TableResults}

\end{table}

We have presented $\icltl$, a variant of Kremer's intuitionistic dynamic topological logic, and shown it to be decidable; we summarize the known decidability results in Table \ref{TableResults}. We also showed that the language of $\icltl$ is expressive enough to characterize minimality and Poincar\'e recurrence, two key properties which sparked interest in $\sf DTL$. This makes $\icltl$ arguably be the first decidable logic suitable for reasoning about non-trivial asymptotic behavior of dynamical topological systems. Nevertheless, our techniques are model-theoretic and do not yield an axiomatization, raising the following:

\begin{question}
Is there a natural axiomatization for $\itlcplus$ and/or its fragments?
\end{question}

Note also that the decision procedure we have given is not elementary, since the size of $\irr\Sigma$ grows super-exponentially on $\lgt\varphi$. Nevertheless, there is little reason to assume that this procedure is optimal; we constructed $\irr \Sigma$ so that a surjective simulation onto any model could always be found, but surjectivity is not necessarily required to falsify a given formula. Hence, a sharp lower bound on the complexity of $\icltl$ remains to be found.

\begin{question}
What is the complexity of the validity problem in $\icltl$?
\end{question}

In case that such a bound is intractable, there are other natural variants of $\sf DTL$ in the spirit of $\icltl$ which may have lower complexity. One could view the topological semantics of intuitionistic logic as a restriction of the classical semantics to the algebra of open sets. However, this is not the only important topologically-defined algebra: the sub-algebra of regular open (or closed) sets has already been studied in the context of spatial reasoning \cite{Duntsch05,Lutz06}, and Lando has studied $\sf DTL$ modulo sets of measure zero \cite{Lando12}. A reduced semantics modulo {\em meager} sets (i.e., countable unions of nowhere-dense sets) would also be meaningful in the context of dynamical systems. For the definitions and basic properties of these algebras, we refer the reader to a text such as~\cite{Givant}.

A separate strategy for finding tractable fragments could be to impose additional syntactical restrictions to $\landi$, such as limiting the number of embedded implications. Minimality is characterized using only one implication, and Poincar\'e recurrence uses three,\footnote{In fact, in view of the intuitionistic tautology $(p\imp \ineg\ineg \nx\ps p)\iiff\ineg(p \wedge \ineg \nx \ps p)$, two implications suffice to characterize Poincar\'e recurrence.} so that such restricted systems might suffice for applications. This strategy has been successfully employed to obtain tractable fragments of the polymodal provability logic $\sf GLP$ \cite{Beklemishev:2013:PositiveProvabilityLogic,Dashkov:2012:PositiveFragment}, which, like $\sf DTL$, is topologically complete but not Kripke complete \cite{BeklemishevGabelaia:2011:TopologicalCompletenessGLP}. This raises the following:

\begin{question}

Can tractable and useful variants of $\itlcplus$ be obtained by

\begin{enumerate}[label=(\alph*)]

\item  using different spatial algebras, or

\item restricting the syntax to suitable fragments?

\end{enumerate}

\end{question}

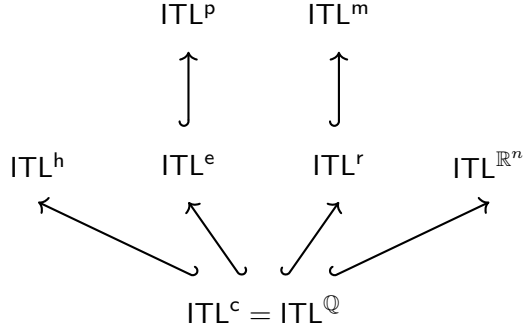
\begin{figure}
\begin{center}

\begin{tikzpicture}

\draw (0,0) node {$\itlcplus= {\sf ITL}^\mathbb Q_{\phantom{\ps}}$};

\draw[left hook->,thick] (-.9,.5) -- (-3,1.5);

\draw[left hook->,thick] (-.3,.5) -- (-1,1.5);

\draw[right hook->,thick] (.3,.5) -- (1,1.5);

\draw[right hook->,thick] (.9,.5) -- (3,1.5);

\draw (-3,2) node {${\sf ITL^h}$};

\draw (-1,2) node {${\sf ITL^e}$};

\draw[right hook->,thick] (-1,2.5) -- (-1,3.5);

\draw (-1,4) node {${\sf ITL^p}$};

\draw (1,2) node {${\sf ITL^r}$};

\draw[right hook->,thick] (1,2.5) -- (1,3.5);

\draw (1,4) node {${\sf ITL^m}$};

\draw (3,2) node {${\sf ITL}^{\mathbb R^n}$};

\end{tikzpicture}

\end{center}

\caption{Inclusions among logics we have considered, with notation as in Table \ref{TableResults}. }
\label{FigCompare}
\end{figure}

Most of the classes of systems we have considered give rise to logics different from $\itlcplus$, as depicted in Figure \ref{FigCompare}. All inclusions in the figure follow from one class being contained in the other, except for ${\sf ITL}^\mathbb Q \subseteq \itlcplus $, which is Theorem \ref{TheoQComplete}. Moreover:
\begin{enumerate}

\item  ${\sf ITL^h}$ is the only logic in the figure containing $\mathord \forall \nx p \imp \forall p$, which is valid over systems with a surjective function.

\item ${\sf ITL^e}$ and ${\sf ITL^p}$ contain $\nec p\imp \nx\nec p$, which is not in ${\sf ITL}^{\mathbb R^n}$ (Example \ref{ExKremer}). They also contain \eqref{EqAleksForm}, which is neither in ${\sf ITL}^{\mathbb R^n}$ nor in ${\sf ITL^h}$, given that the function in Example \ref{ExQ} uses a homeomorphism.

\item ${\sf ITL^h}$ and ${\sf ITL^p}$ are the only logics containing $(\nx p \imp \nx q) \imp \nx (p\imp q)$; this formula is valid over systems with an open function.

\item ${\sf ITL^r}$ and ${\sf ITL^m}$ are the only logics containing \eqref{EqPoincInt} (Theorem \ref{TheoPoincare}),

\item ${\sf ITL^m}$ is the only logic containing \eqref{EqMinFor} (Theorem \ref{TheoITLMinExpress}), and

\item ${\sf ITL}^{\mathbb R^n}$ is the only logic containing \eqref{EqConn} \cite{ShehtmanEverywhereHere}.

\end{enumerate}
We leave a full proof of these claims to the reader; one useful counter-example for this is given by Figure \ref{FigMinNonInv}, which presents a disconnected, non-invertible Alexandroff minimal system. Only the following questions are left open by this analysis:
\begin{question}\label{QuestInclus} Do the inclusions ${\sf ITL^e} \subseteq {\sf ITL^r} $ or ${\sf ITL^e} \subseteq {\sf ITL^m} $ hold?
\end{question}

Other than Question \ref{QuestInclus}, Figure \ref{FigCompare} is complete (in the sense that all inclusions are shown), and remains complete if we replace the logics by the respective $\landi$ fragments.
Given that these logics are mostly distinct, it is an interesting open problem whether intuitionistic temporal logics over special classes of systems are more feasible than their classical counterparts. For example, $\DTLm$ is decidable, unlike the unrestricted $\sf DTL$: perhaps $\ITLm_\ps$ also has lower complexity than $\icltl$? Similarly, it is not known if $\sf DTL$ over Poincar\'e recurrent systems is decidable, but settling the decidability of intuitionistic temporal logic over this class may be a more accessible problem.

\begin{question}
Which of the unknown logics of Table \ref{TableResults} are decidable, and what is their complexity?
\end{question}

Finally, we leave open the question of the decidability of Kremer's original logic over the $\{\circ,\nec\}$-fragment.

\begin{question}
Is $\itlk$ decidable?
\end{question}

\section*{Acknowledgments}

The author would like to thank Philip Kremer for permission to reproduce Example \ref{ExKremer}, as well as the anonymous referee for his or her helpful comments.

 

\begin{thebibliography}{10}

\bibitem{akin}
E.~Akin.
\newblock {\em The General Topology of Dynamical Systems}.
\newblock Graduate Studies in Mathematics. American Mathematical Society, 1993.

\bibitem{aleksandroff}
P.~Alexandroff.
\newblock Diskrete {R}\"{a}ume.
\newblock {\em Matematicheskii Sbornik}, 2(44):501--518, 1937.

\bibitem{arte}
S.~Artemov, J.~Davoren, and A.~Nerode.
\newblock Modal logics and topological semantics for hybrid systems.
\newblock {\em Technical Report MSI 97-05}, 1997.

\bibitem{Balbiani2017}
P.~{Balbiani}, J.~{Boudou}, M.~{Di{\'e}guez}, and D.~{Fern{\'a}ndez-Duque}.
\newblock Bisimulations for intuitionistic temporal logics.
\newblock {\em arXiv}, 1803.05078, 2018.

\bibitem{BalbianiDieguezJelia}
P.~Balbiani and M.~Di\'eguez.
\newblock Temporal here and there.
\newblock In M.~Loizos and A.~Kakas, editors, {\em Logics in Artificial
  Intelligence}, pages 81--96. Springer, 2016.

\bibitem{Balibrea}
F.~Balibrea, T.~Downarowicz, R.~Hric, L.~Snoha, and V.~\v{S}pitalsk\'y.
\newblock Almost totally disconnected minimal systems.
\newblock {\em Ergodic Theory and Dynamical Systems}, 29(3):737--766, 2009.

\bibitem{Beklemishev:2013:PositiveProvabilityLogic}
L.~Beklemishev.
\newblock Positive provability logic for uniform reflection principles.
\newblock {\em Annals of Pure and Applied Logic}, 165(1):82--105, 2014.

\bibitem{BeklemishevGabelaia:2011:TopologicalCompletenessGLP}
L.~Beklemishev and D.~Gabelaia.
\newblock Topological completeness of the provability logic $\mathsf{GLP}$.
\newblock {\em Annals of Pure and Applied Logic}, 164(12):1201--1223, 2013.

\bibitem{Bir}
G.~Birkhoff.
\newblock Quelques th\'eor\`emes sur le mouvement des syst\`emes dynamiques.
\newblock {\em Bulletin de la Soci\'et\'e Math\'ematique de France},
  40:305--323, 1912.

\bibitem{Boudou2017}
J.~Boudou, M.~Di\'eguez, and D.~Fern\'andez-Duque.
\newblock A decidable intuitionistic temporal logic.
\newblock In {\em 26th {EACSL} Annual Conference on Computer Science Logic,
  {CSL} 2017, August 20-24, 2017, Stockholm, Sweden}, pages 14:1--14:17, 2017.

\bibitem{Bruin}
H.~Bruin, S.~Kolyada, and L'. Snoha.
\newblock Minimal nonhomogeneous continua.
\newblock {\em Colloquium Mathematicum}, 95(1):123--132, 2003.

\bibitem{Dashkov:2012:PositiveFragment}
E.~Dashkov.
\newblock On the positive fragment of the polymodal provability logic
  $\mathsf{GLP}$.
\newblock {\em Mathematical Notes}, 91(3-4):318--333, 2012.

\bibitem{Davies96}
R.~Davies.
\newblock A temporal-logic approach to binding-time analysis.
\newblock In {\em Proceedings, 11th Annual {IEEE} Symposium on Logic in
  Computer Science, New Brunswick, New Jersey, USA, July 27-30, 1996}, pages
  184--195, 1996.

\bibitem{Davoren09}
J.~Davoren.
\newblock On intuitionistic modal and tense logics and their classical
  companion logics: Topological semantics and bisimulations.
\newblock {\em Annals of Pure and Applied Logic}, 161(3):349--367, 2009.

\bibitem{DavorenIntuitionistic}
J.~Davoren, V.~Coulthard, T.~Moor, R.~Gor\'e, and A.~Nerode.
\newblock Topological semantics for intuitionistic modal logics, and spatial
  discretisation by {A}/{D} maps.
\newblock In {\em Workshop on Intuitionistic Modal Logic and Applications
  (IMLA)}, 2002.

\bibitem{Duntsch05}
I.~D{\"{u}}ntsch and M.~Winter.
\newblock A representation theorem for {B}oolean contact algebras.
\newblock {\em Theoretical Computer Science}, 347(3):498--512, 2005.

\bibitem{Ewald}
W.~Ewald.
\newblock Intuitionistic tense and modal logic.
\newblock {\em The Journal of Symbolic Logic}, 51(1):166--179, 1986.

\bibitem{me2}
D.~Fern{\'a}ndez-Duque.
\newblock Non-deterministic semantics for dynamic topological logic.
\newblock {\em Annals of Pure and Applied Logic}, 157(2-3):110--121, 2009.

\bibitem{me:metric}
D.~Fern{\'a}ndez-Duque.
\newblock Dynamic topological logic interpreted over metric spaces.
\newblock {\em Journal of Symbolic Logic}, 2011.

\bibitem{me:minimal}
D.~Fern{\'a}ndez-Duque.
\newblock Dynamic topological logic interpreted over minimal systems.
\newblock {\em Journal of Philosophical Logic}, 40(6):767--804, 2011.

\bibitem{dtlaxiom}
D.~Fern\'{a}ndez-Duque.
\newblock A sound and complete axiomatization for dynamic topological logic.
\newblock {\em Journal of Symbolic Logic}, 77(3):947--969, 2012.

\bibitem{DFDTOCL}
D.~Fern{\'{a}}ndez{-}Duque.
\newblock Non-finite axiomatizability of dynamic topological logic.
\newblock {\em {ACM} Transactions on Computational Logic}, 15(1):4:1--4:18,
  2014.

\bibitem{pml}
D.~Gabelaia, A.~Kurucz, F.~Wolter, and M.~Zakharyaschev.
\newblock Non-primitive recursive decidability of products of modal logics with
  expanding domains.
\newblock {\em Annals of Pure and Applied Logic}, 142(1-3):245--268, 2006.

\bibitem{Givant}
S.~Givant and P.~Halmos.
\newblock {\em Introduction to Boolean Algebras}.
\newblock Undergraduate Texts in Mathematics. Springer, New York, 2009.

\bibitem{Godel1931}
K.~G\"{o}del.
\newblock {\"{U}ber Formal Unentscheidbare S\"{a}tze der Principia Mathematica
  und Verwandter Systeme, I}.
\newblock {\em Monatshefte f\"{u}r Mathematik und Physik}, 38:173--198, 1931.

\bibitem{KamideBounded}
N.~Kamide and H.~Wansing.
\newblock Combining linear-time temporal logic with constructiveness and
  paraconsistency.
\newblock {\em Journal of Applied Logic}, 8(1):33--61, 2010.

\bibitem{KojimaNext}
K.~Kojima and A.~Igarashi.
\newblock Constructive linear-time temporal logic: Proof systems and {K}ripke
  semantics.
\newblock {\em Information and Computation}, 209(12):1491 -- 1503, 2011.

\bibitem{konev}
B.~Konev, R.~Kontchakov, F.~Wolter, and M.~Zakharyaschev.
\newblock Dynamic topological logics over spaces with continuous functions.
\newblock In G.~Governatori, I.~Hodkinson, and Y.~Venema, editors, {\em
  Advances in Modal Logic}, volume~6, pages 299--318, London, 2006. College
  Publications.

\bibitem{wolter}
B.~Konev, R.~Kontchakov, F.~Wolter, and M.~Zakharyaschev.
\newblock On dynamic topological and metric logics.
\newblock {\em Studia Logica}, 84:129--160, 2006.

\bibitem{KremerIntuitionistic}
P.~Kremer.
\newblock A small counterexample in intuitionistic dynamic topological logic,
  2004.
\newblock
  http://individual.utoronto.ca/philipkremer/onlinepapers/counterex.pdf.

\bibitem{s5}
P.~Kremer.
\newblock Dynamic topological $\mathsf{S5}$.
\newblock {\em Annals of Pure and Applied Logic}, 160:96--116, 2009.

\bibitem{kmints}
P.~Kremer and G.~Mints.
\newblock Dynamic topological logic.
\newblock {\em Annals of Pure and Applied Logic}, 131:133--158, 2005.

\bibitem{mdml}
A.~Kurucz, F.~Wolter, M.~Zakharyaschev, and D.~Gabbay.
\newblock {\em {Many-Dimensional Modal Logics: Theory and Applications}},
  volume 148 of {\em Studies in Logic and the Foundations of Mathematics}.
\newblock North Holland, first edition, 2003.

\bibitem{Lando12}
T.~Lando.
\newblock Dynamic measure logic.
\newblock {\em Annals of Pure and Applied Logic}, 163(12):1719--1737, 2012.

\bibitem{temporal}
O.~Lichtenstein and A.~Pnueli.
\newblock Propositional temporal logics: Decidability and completeness.
\newblock {\em Logic Jounal of the IGPL}, 8(1):55--85, 2000.

\bibitem{Lutz06}
C.~Lutz and F.~Wolter.
\newblock Modal logics of topological relations.
\newblock {\em Logical Methods in Computer Science}, 2(2), 2006.

\bibitem{Munkres}
J.~Munkres.
\newblock {\em Topology}.
\newblock Prentice Hall, second edition, 2000.

\bibitem{NishimuraConstructivePDL}
H.~Nishimura.
\newblock Semantical analysis of constructive {P}{D}{L}.
\newblock {\em Publications of the Research Institute for Mathematical
  Sciences, Kyoto University}, 18:427--438, 1982.

\bibitem{Poincare}
H.~Poincar\'e.
\newblock Sur le probl\`eme des trois corps et les \'equations de la dynamique.
\newblock {\em Acta Mathematica}, 13:1--270, 1890.

\bibitem{ShehtmanEverywhereHere}
V.~Shehtman.
\newblock `{E}verywhere' and `here'.
\newblock {\em Journal of Applied Non-Classical Logics}, 9(2-3):369--379, 1999.

\bibitem{Simpson94}
A.~Simpson.
\newblock {\em The proof theory and semantics of intuitionistic modal logic}.
\newblock PhD thesis, University of Edinburgh, {UK}, 1994.

\bibitem{tarski}
A.~Tarski.
\newblock Der {A}ussagenkalk\"ul und die {T}opologie.
\newblock {\em Fundamenta Mathematica}, 31:103--134, 1938.

\end{thebibliography}

\end{document}